\documentstyle[12pt]{article}
\makeatletter
\def\emline#1#2#3#4#5#6{%
       \put(#1,#2){\special{em:moveto}}%
       \put(#4,#5){\special{em:lineto}}}
\def\hybrid{\topmargin 0pt      \oddsidemargin 0pt
        \headheight 0pt \headsep 0pt
        \textwidth 6.35in       % A4 paper
        \textheight 9.5in       % A4 paper
        \marginparwidth 0.0in
        \parskip 5pt plus 1pt   \jot = 1.5ex}
\catcode`\@=11
\def\marginnote#1{}

\newcount\hour
\newcount\minute
\newtoks\amorpm
\hour=\time\divide\hour by60
\minute=\time{\multiply\hour by60 \global\advance\minute by-\hour}
\edef\standardtime{{\ifnum\hour<12 \global\amorpm={am}%
        \else\global\amorpm={pm}\advance\hour by-12 \fi
        \ifnum\hour=0 \hour=12 \fi
        \number\hour:\ifnum\minute<10 0\fi\number\minute\the\amorpm}}
\edef\militarytime{\number\hour:\ifnum\minute<10 0\fi\number\minute}

\def\draftlabel#1{{\@bsphack\if@filesw {\let\thepage\relax
   \xdef\@gtempa{\write\@auxout{\string
      \newlabel{#1}{{\@currentlabel}{\thepage}}}}}\@gtempa
   \if@nobreak \ifvmode\nobreak\fi\fi\fi\@esphack}
        \gdef\@eqnlabel{#1}}
\def\@eqnlabel{}
\def\@vacuum{}
\def\draftmarginnote#1{\marginpar{\raggedright\scriptsize\tt#1}}

\def\draftlabel#1{{\@bsphack\if@filesw {\let\thepage\relax
   \xdef\@gtempa{\write\@auxout{\string
      \newlabel{#1}{{\@currentlabel}{\thepage}}}}}\@gtempa
   \if@nobreak \ifvmode\nobreak\fi\fi\fi\@esphack}
        \gdef\@eqnlabel{#1}}
\def\@eqnlabel{}
\def\@vacuum{}
\def\draftmarginnote#1{\marginpar{\raggedright\scriptsize\tt#1}}

\def\draft{\oddsidemargin -.5truein
        \def\@oddfoot{\sl preliminary draft \hfil
        \rm\thepage\hfil\sl\today\quad\militarytime}
        \let\@evenfoot\@oddfoot \overfullrule 3pt
        \let\label=\draftlabel
        \let\marginnote=\draftmarginnote
   \def\@eqnnum{(\theequation)\rlap{\kern\marginparsep\tt\@eqnlabel}%
\global\let\@eqnlabel\@vacuum}  }

\def\numberbysection{\@addtoreset{equation}{section}
        \def\theequation{\thesection.\arabic{equation}}}

\def\underline#1{\relax\ifmmode\@@underline#1\else
        $\@@underline{\hbox{#1}}$\relax\fi}

\def\titlepage{\@restonecolfalse\if@twocolumn\@restonecoltrue\onecolumn
     \else \newpage \fi \thispagestyle{empty}\c@page\z@
        \def\thefootnote{\fnsymbol{footnote}} }

\def\endtitlepage{\if@restonecol\twocolumn \else  \fi
        \def\thefootnote{\arabic{footnote}}
        \setcounter{footnote}{0}}  %\c@footnote\z@ }
\relax

\numberbysection
\hybrid

\newfont{\Bbb}{msbm10 scaled 1\@ptsize00}
\newcommand{\CC}{\mbox{\Bbb C}}

\newcommand{\RR}{\mbox{\Bbb R}}

\newcommand{\ZZ}{\mbox{\Bbb Z}}
\newfont{\Bbbb}{msbm7 scaled 1\@ptsize00}
\newcommand{\z}{\raise-1pt\hbox{$\mbox{\Bbbb Z}$}}
\def\beq{\begin{equation}}
\def\eeq{\end{equation}}
\def\p{\partial}

\def\normord{ {\scriptstyle {{\bullet}\atop{\bullet}}} }

\begin{document}

\begin{titlepage}

\title{Commuting difference operators with elliptic
coefficients from Baxter's vacuum vectors}

\author{A. Zabrodin
\thanks{Joint Institute of Chemical Physics, Kosygina str. 4, 117334,
Moscow, Russia and ITEP, 117259, Moscow, Russia}}
\date{December 1999}

\maketitle

\begin{abstract}

For quantum integrable models
with elliptic $R$-matrix, we
construct the Baxter $Q$-operator in infinite-dimensional
representations of the algebra of observables.

\end{abstract}

%\vfill

\end{titlepage}

\section{Introduction}

The method of $Q$-operator, first proposed by Baxter
in his solution of the
8-vertex model and the spin-$\frac{1}{2}$ $XYZ$ chain \cite{Baxter},
is perhaps the most powerful approach to quantum integrable
models known today. As compared with the
algebraic Bethe ansatz \cite{FT},
one or another version of this approach,
known also as the method of vacuum vectors,
is indispensable for models which do not
possess bare vacuum states.
Among later striking developments we mention the works \cite{CP},
where the Baxter method was applied
to the chiral Potts model.

At that time, however, the $Q$-operator itself
seemed to be a kind of virtual
object for its explicit form was not known.
In Baxter's original work, the $Q$-operator played an auxiliary
role since the local Hamiltonian
to be diagonalized was made of
the quantum transfer matrix ${\bf T}(\lambda )$.
The main properties of the $Q$-operator
${\bf Q}(\lambda )$ are the commutativity
$[{\bf Q}(\lambda ) ,\,{\bf Q}(\mu )]=0$,
$[{\bf Q}(\lambda ) ,\,{\bf T}(\mu )]=0$ for all
$\lambda$, $\mu$, and the following $T$-$Q$-relation:
\beq
\label{01}
{\bf T}(\lambda ){\bf Q}(\lambda )=
a(\lambda ){\bf Q}(\lambda +\eta)+b(\lambda )
{\bf Q}(\lambda -\eta)\,.
\eeq
Here $a(\lambda )$, $b(\lambda )$ are $c$-number functions
depending on a particular model. Combined with
(model-dependent) analytical
properties of the ${\bf Q}(\lambda )$, this relation
allows one to obtain a system of Bethe equations and
eventually to solve the diagonalization problem.
In the absense of global bare vacuum states, Baxter's method
to derive (\ref{01}) was to consider a
continuous-parametric family of
special states that factorize into products of
{\it local} vacuum vectors. Such states have especially
simple properties under action of ${\bf T}(\lambda )$.
Treating them as columns or rows of some matrix
${\bf Q}(\lambda )$, one arrives, after a proper normalization,
at (\ref{01}). Some years ago this approach was extended \cite{YB}
to the case of twisted boundary conditions.
It remained obscure, however, whether the $Q$-matrix could play any
other role than just a technical tool for deriving Bethe equations.

The situation changed only quite recently, when it was
realized that
in integrable lattice models with an infinite-dimensional
Hilbert space at each site it is ${\bf Q}$ rather than
${\bf T}$ that is the
generating function of local Hamiltonians.
With this motivation in mind, a considerable
progress in constructing the $Q$-operator for a number
of integrable models and in understanding
its nature was achieved.

The work \cite{pg}
was perhaps the first, where the Baxter approach
was adopted to integrable models with an infinite-dimensional
quantum space. Specifically, in that paper
the $Q$-operator for the quantum Toda chain was
constructed as
certain integral operator with explicitly known kernel.
Recently, this remarkable development
was extended \cite{derk}
to the integrable homogeneous $XXX$ spin chain in a unitary
discrete series representation
of the non-compact group $SU(1,1)$ (see also \cite{pronko}).

For integrable models of $XXZ$-type, the nature of
$Q$-operator was clarified in \cite{BLZ},\,\cite{AF},
where it was given a representation theory meaning.
Specifically, the $Q$-operator can be
obtained from the universal $R$-matrix
for the quantum affine algebra $U_q(\widehat{sl(2)})$,
by taking trace of a chain product of the universal
$R$-matrices in a special $q$-oscillator representation
in one of the spaces.
The commutativity of the $Q$-operators, which comes
out as a surprize in Baxter's approach, is now guaranteed
by characteristic properties of the universal $R$-matrix.
The $T$-$Q$ relation, from this point of
view, is a simple consequence of
fusion rules for representations of the quantum
affine algebra. In fact, two different commuting
$Q$-operators, ${\bf Q}^{+}(\lambda )$ and
${\bf Q}^{-}(\lambda )$, can be introduced this way.
They are two linearly independent operator solutions to the
$T$-$Q$-relation (\ref{01}) treated as a second order
difference equation in $\lambda$. As such, they obey
a quantum Wronskian condition, which turns out to be
even more instructive than the $T$-$Q$ relation itself.

As for matters more customary
for theory of solitons,
it was argued in \cite{ks,kss} that
$Q$-operator is a quantum analog of B\"acklund
transformation. In other words,
it defines an integrable quantum model
with discrete time by
setting ${\bf Q}(\lambda )$ to be the time shift operator.
For any observable ${\bf O}$, one sets
${\bf O}(t+\delta )={\bf Q}(\lambda ){\bf O}(t)
({\bf Q}(\lambda ))^{-1}$, where the discrete time
shift $\delta$ is a function of $\lambda$. Integrability
implies that these discrete evolution operators commute
for all $\delta$.

However, all this was about models
associated with rational or trigonometric
solutions of the Yang-Baxter equation ($R$-matrices).
General $R$-matrices
with additive spectral parameter are
parametrized via elliptic functions.
The symplest elliptic $R$-matrix is
\beq
\label{02}
{\sf R}(\lambda )=\sum_{a=0}^{3}
\frac{\theta_{a+1}(2\lambda +\eta |\tau )}
{\theta_{a+1}(\eta |\tau )}
\sigma_a \otimes \sigma_a \,,
\eeq
where $\theta_a(z|\tau)$ are Jacobi $\theta$-functions,
$\sigma_a$ are Pauli matrices, and $\sigma_0$ is the unit
matrix.
This $R$-matrix is associated with the 8-vertex model, being the
matrix of local Boltzmann weights at a vertex.
The transfer matrix of this model is the generating function
of conserved quantities for
the integrable $XYZ$ spin-$\frac{1}{2}$ chain.
By models of the $XYZ$-type we understand more general
integrable models with elliptic $R$-matrix, for instance
integrable anisotropic spin chains with higher spin.
Integrable spin chains of $XYZ$-type can be solved
by the generalized algebraic Bethe
ansatz \cite{FT},\,\cite{Takebe}.

However, the quantum affine algebra approach is not applicable
in the $XYZ$-case since
nothing similar to the universal elliptic $R$-matrix is known.
Nevertheless, the elliptic dependence on the spectral parameter
does not destroy the functional relations
between transfer matrices and $Q$-operators,
including the $T$-$Q$ relation
and the Wronskian condition \cite{BR},\,\cite{KLWZ}.
This suggests that everything what is possible to do
in $XXX$ and $XXZ$ cases, might be possible,
mutatis mutandis, for $XYZ$, too.

In this paper we present a direct construction of the
$Q$-operator for homogeneous $N$-site spin chains of the $XYZ$-type.
In lattice integrable models with elliptic $R$-matrix,
the algebra of local observables
is the Sklyanin algebra \cite{Skl1,Skl2}. A concrete
model is defined by fixing a particular representation
of this algebra. Similar to the $sl(2)$-case
(models of $XXX$-type), the representations are labeled
by a continuous parameter which is called
{\it spin}, by analogy with $sl(2)$. We allow the
spin to take any complex value, so we
work in an infinite-dimensional representation of the
algebra of observables. In other words, we deal with
representations of what would be called
a non-compact real form of the Sklyanin algebra.
Such representations can be realized by difference
operators, so our $Q$-operator is a diffeference
operator in many variables.

Our strategy closely follows
the original Baxter argument.
First, we construct a family of local
(infinite-dimensional) vacuum vectors that obey
the celebrated ``pair propagation through
a vertex" property. Next, their product over all sites
of the chain is, up to normalization, the kernel
of the $Q$-operator. At last, one should prove
the commutativity of the $Q$-operators.
However, the practical realization of these ideas in the
infinite-dimensional setting is by no means obvious.
Technically, it is very different
from what is customary in the 8-vertex model.
The construction extensively uses such truly special functions as
elliptic generalizations of the
gamma-function and basic hypergeometric functions,
and highly non-trivial identities for them.
For example, commutativity of the $Q$-operators
relies on an elliptic version
of the Bailey transformation
formula for the very-well-poised basic hypergeometric
series ${}_{10}W_{9}$. The appearence of
such a sophisticated mathematics in integrable models
is not so surprizing.
Recently, the elliptic series ${}_{10}W_{9}$
was shown \cite{FrTur} to offer a natural representation for
``elliptic $6$-$j$-symbols", which are nothing else than
solutions of the Yang-Baxter equation in a special basis.

It is now a challenge
to reveal a representation theory meaning
of the elliptic $Q$-operator. Hopefully, this would make
its properties and their proofs obvious.

The paper is organized as follows. In the next section
we present the main result in the shortest possible form.
Sec.\,3 contains the necessary material on rerpesentations
of the Sklyanin algebra. The construction of the $Q$-operator
begins in Sec.\,4, where we find a family of
local vacuum states for the quantum ${\sf L}$-operator.
In Sec.\,5 we investigate direct products of such states
along the chain and define unnormalized
``pre-$Q$-operators", ${\bf Q}_{L}$ and ${\bf Q}_{R}$, which
do not yet commute.
Sec.\,6 is rather technical. It is devoted to the
proof of the basic commutation relation for the pre-$Q$-operators,
${\bf Q}_{L}(\lambda ){\bf Q}_{R}(\mu )=
{\bf Q}_{L}(\mu ){\bf Q}_{R}(\lambda )$,
which entails commutativity of the $Q$-operators
introduced in Sec.\,7.
In Sec.\,8 we discuss more specific issues: a) twisted
boundary conditions, b) the
exceptional case of a ``spin chain" with only one site, $N=1$, and
c) the trigonometric ($XXZ$-type) degeneration.
A few concluding remarks are given in Sec.\,9.

Appendix A contains necessary information on the special
functions involved in the construction. In Appendix B we give
a direct proof of the intertwining relation between
representations of spins $\ell$ and $-\ell -1$.

\section{The main result}

To state the result, we recall the definition of the
quantum transfer matrix ${\bf T}(\lambda )$ for the
$XYZ$ spin chain with spin $\ell$.
The elliptic quantum $\mbox{{\sf L}}$-operator is
the matrix
\beq
\mbox{{\sf L}}(\lambda)=\frac{1}{2}
\left ( \begin{array}{cc}
\theta _{1}(2\lambda){\bf s}_0 +
\theta _{4}(2\lambda){\bf s}_3 &
\theta _{2}(2\lambda){\bf s}_1 +
\theta _{3}(2\lambda){\bf s}_2
\\& \\
\theta _{2}(2\lambda){\bf s}_1 -
\theta _{3}(2\lambda){\bf s}_2 &
\theta _{1}(2\lambda){\bf s}_0 -
\theta _{4}(2\lambda){\bf s}_3
\end{array} \right )
\label{L}
\eeq
with non-commutative matrix elements. Specifically,
${\bf s}_a$ are difference operators in a
complex variable $z$:
\beq
{\bf s}_{a} = \frac{\theta _{a+1}(2z -2\ell \eta)}
{\theta _{1}(2z)}\, e^{\eta \p _z}
-\frac{\theta _{a+1}(-2z -2\ell \eta)}
{\theta _{1}(2z)}\, e^{-\eta \p _z}\,,
\label{Sa}
\eeq
introduced by Sklyanin \cite{Skl2}.
Here $\theta_a(z)\equiv \theta_a(z|\tau)$ are
Jacobi $\theta$-functions
with the elliptic module $\tau$,
$\mbox{Im}\,\tau >0$, $\ell$ is a complex number (the spin),
and $\eta \in \CC$ is a parameter which is assumed to
belong to the fundamental parallelogram
with vertices $0$, $1$, $\tau$, $1+\tau$, and to be
incommensurate with $1, \tau$.
Definitions and transformation properties
of the $\theta$-functions are listed in Appendix A.

Let us consider the matrix product of
${\sf L}$-operators (\ref{L}):
$\stackrel{1}{{\sf L}}(\lambda)
\stackrel{2}{{\sf L}}(\lambda) \ldots
\stackrel{N}{{\sf L}}(\lambda)$, $N>1$,
where each operator matrix
$\stackrel{i}{{\sf L}}(\lambda)$ has the form (\ref{L})
with ${\bf s}_a$ acting on its own variable $z_i$.
Clearly, matrix elements of
$\stackrel{i}{{\sf L}}(\lambda)$ and
$\stackrel{j}{{\sf L}}(\lambda)$ commute for $i\neq j$.
The trace of this product (the sum of two diagonal elements),
\beq
\label{T3}
{\bf T}(\lambda )=\mbox{tr}\, \bigl [
\stackrel{1}{{\sf L}}(\lambda)
\stackrel{2}{{\sf L}}(\lambda) \, \ldots \,
\stackrel{N}{{\sf L}}(\lambda) \bigr ]\,,
\eeq
is called the quantum transfer matrix of the spin chain.

Our main result is the explicit realization of
the $Q$-operator for this spin chain:
\beq
\label{QQ}
{\bf Q}(\lambda )=
\Gamma^{-N}(2\lambda -2\ell \eta )
\, \normord
\prod_{i=1}^{N}{\bf q}_{i, i\!+\!1}(\lambda )
\normord \,
e^{(\lambda +\ell \eta )\sum_{j=1}^{N}\p _{z_j}}\,,
\eeq
where
\beq
\label{qq}
\begin{array}{lll}
{\bf q}_{i, i\!+\!1}(\lambda )
&=&\!\!
\displaystyle{
\frac{\Gamma (2z_i )\,\Gamma (z_i \!-\!z_{i+1}\!+\!2\lambda )
\Gamma (z_i \!+\!z_{i+1}\!+\!2\lambda )}
{\Gamma (2z_i \!+\!2\lambda \!+\!2\ell \eta )
\Gamma (z_i \!-\!z_{i+1}\!-\!2\ell \eta )
\Gamma (z_i \!+\!z_{i+1}\!-\!2\ell \eta )}}
\,\, \times
\\ && \\
&\times & \!\! \!
\displaystyle{
{}_{6}W_{5} \Bigl (
-\frac{z_i\!+\!\lambda }{\eta}\!-\!\ell ;\,
-\frac{\lambda }{\eta}\!-\!\ell ,\,
-\frac{z_{i+1}\!+\!z_i}{2\eta}\!-\!\ell ,\,
\frac{z_{i+1}\!-\!z_i}{2\eta}\!-\!\ell ; \,
e^{-2\eta \p _{z_i}} }\Bigl | \eta , \tau \Bigr )
\end{array}
\eeq
(it is implied that $z_{N+1}=z_1$).
Here $\Gamma (z)=\Gamma (z|\tau , 2\eta )$ is the elliptic
gamma-function \cite{R3,FV3}:
$$
\Gamma (z|\tau , 2\eta )=\! \prod_{k,m =0}^{\infty}
\frac{1\!-\!e^{2\pi i ( (k+1)\tau +2(m+1)\eta -z)}}
{1\!-\!e^{2\pi i ( k\tau +2m\eta +z)}}
\!=\! \exp \left (
-i\sum_{k=1}^{\infty} \frac{\sin (\pi k (2z\!-\!2\eta \!-\!\tau ))}
{2k\sin (2\pi k \eta )\sin (\pi k \tau )} \right ),
$$
and ${}_{6}W_{5}(\alpha_1; \alpha_4,
\ldots , \alpha_6; \,z\, |\eta , \tau )$ is the elliptic
analog of the very-well-poised
hypergeometric series (with the operator argument):
$$
{}_{6}W_{5}(\alpha_1; \alpha_4,
\ldots , \alpha_6; \,z\, |\eta , \tau )=
\sum_{k=0}^{\infty}
\frac{[\alpha_1 +2k ][\alpha_1 ]_k}
{[\alpha_1][1]_k} \Bigl (
\prod_{m=1}^{3}\frac{[\alpha_{m+3}]_k}
{[\alpha_1 \!-\! \alpha_{m+3} \!+\! 1]_k} \Bigr )\, z^k\,,
$$
where $[\alpha ]_k =\prod_{j=0}^{k-1}
\theta_1 (2\alpha \eta \!+\!2j\eta )$.
For more information on the special functions
$\Gamma$, ${}_{6}W_{5}$ see Appendix A.
The normal ordering means that all operators
$e^{-2\eta \p_{z_j}}$ are moved to the right.

We prove the commutativity
$[{\bf Q}(\lambda ),\, {\bf T}(\mu )]=[{\bf Q}(\lambda ),
\, {\bf Q}(\mu )]=0$, and the $T$-$Q$ relation
\beq
\label{03}
{\bf T}(\lambda ){\bf Q}(\lambda )=
\theta_{1}^{N}(2\lambda -2\ell \eta )
{\bf Q}(\lambda +\eta )
+\theta_{1}^{N}(2\lambda +2\ell \eta )
{\bf Q}(\lambda -\eta )\,.
\eeq
Another solution to this operator difference equation has
a similar form (see (\ref{N7})).

The limiting case $\tau \to 0$, when elliptic
functions degenerate to trigonometric ones, corresponds to the
$XXZ$ spin chain.
The structure of the formulas (\ref{QQ}), (\ref{qq})
remains the same while the special functions get simplified.
The elliptic gamma-function
becomes the $q$-gamma function (or the $q$-exponential function)
while the elliptic hypergeometric series ${}_{6}W_{5}$
turns into the little $q$-Jacobi function represented by
the basic hypergeometric series ${}_{2}\phi _{1}$.

\section{Representations of the Sklyanin algebra}

The aim of this section is to give the necessary
preliminaries on representations of the Sklyanin algebra.
We begin with a few formulas related to
the quantum ${\sf L}$-operator
with elliptic dependence on the spectral parameter.

Consider the quantum ${\sf L}$-operator (\ref{L}),
where ${\bf s}_a$ are realized as in (\ref{Sa}).
The four operators ${\bf s}_a$ obey the commutation
relations of the Sklyanin algebra\footnote{The standard generators
of the Sklyanin algebra \cite{Skl1}, $S_a$,
are related to ours as follows:
$S_{a}=(i)^{\delta _{a,2}}\theta _{a+1}(\eta ){\bf s}_a$.}:
\beq
\begin{array}{l}
(-1)^{\alpha +1}I_{\alpha 0}{\bf s}_{\alpha}{\bf s}_{0}=
I_{\beta \gamma}{\bf s}_{\beta}{\bf s}_{\gamma}
-I_{\gamma \beta}{\bf s}_{\gamma}{\bf s}_{\beta}\,,
\\ \\
(-1)^{\alpha +1}I_{\alpha 0}{\bf s}_0 {\bf s}_{\alpha}=
I_{\gamma \beta}{\bf s}_{\beta}{\bf s}_{\gamma}
-I_{\beta \gamma}{\bf s}_{\gamma}{\bf s}_{\beta}
\end{array}
\label{skl6}
\eeq
with the structure constants
$I_{ab}=\theta _{a+1}(0)\theta _{b+1}(2\eta)$.
Here $a,b =0, \ldots , 3$ and
$\{\alpha ,\beta , \gamma \}$ stands for any cyclic
permutation of $\{1 ,2,3\}$.
The relations of the Sklyanin algebra
are equivalent to the condition that the $\mbox{{\sf L}}$-operator
satisfies the
``$\mbox{{\sf R}}\mbox{{\sf L}}
\mbox{{\sf L}}=\mbox{{\sf L}}
\mbox{{\sf L}}\mbox{{\sf R}}$"
relation with the elliptic $R$-matrix (\ref{02}).

Plugging the difference operators (\ref{Sa}) into (\ref{L}),
one can represent the ${\sf L}$-operator in a factorized form,
which is especially convenient for computations:
\beq
\label{fact}
{\sf L}(\lambda)=\theta_{1}(2\lambda_{+})
\normord
V^{-1}(\lambda_{+},z)
\left ( \begin{array}{cc}
e^{\eta \p_{z}}&0 \\
0 & e^{-\eta \p_{z}} \end{array} \right )
V(\lambda_{-},z)\normord \,.
\eeq
Hereafter we use the convenient notation
$\lambda_{\pm}\equiv \lambda \pm \ell \eta$,
$\bar \theta_{a}(z)\equiv \theta_{a}(z|\frac{\tau}{2})$.
The dots mean the normal ordering when the shift operators
are moved to the right and
$$
V(\lambda , z)=\left ( \begin{array}{cc}
\bar \theta_{4}(z+\lambda )&
\bar \theta_{3}(z+\lambda ) \\
\bar \theta_{4}(z-\lambda )&
\bar \theta_{3}(z-\lambda ) \end{array} \right ).
$$
Note that
\beq
\label{trans}
{\sf L}(-\lambda +{\scriptstyle \frac{1}{2} }
(1 +\tau))=-
e^{-\pi i \tau +4\pi i \lambda}
{\sf L}^{\top}(\lambda )
\eeq
where
${\sf L}^{\top}$ means transposition in the two-dimensional
space.

The parameter
$\ell$ in (\ref{Sa}) is called the spin of the representation.
If necessary, we write ${\bf s}_a ={\bf s}_{a}^{(\ell )}$
to indicate the dependence on $\ell$.
When $\ell \in
\frac{1}{2} \ZZ_{+}$, these operators
have a finite-dimensional invariant subspace, namely,
the space
$\Theta_{4\ell}^{+}$ of {\it even} $\theta$-functions of
order $4\ell$ (see Appendix A). This
is the representation space of the $(2\ell +1)$-dimensional
irreducible representation (of series a))
of the Sklyanin algebra. For example,
at $\ell =\frac{1}{2}$ the functions
$\bar \theta_4 (z)$, $\bar \theta_3 (z)$ form a basis
in $\Theta_{2}^{+}$, and the generators
${\bf s}_a$, with respect to this basis,
are represented by $2$$\times$$2$
matrices $(-i)^{\delta_{a,2}}
(\theta_{a+1}(\eta ))^{-1}\sigma_a$. In this case,
${\sf L}(\lambda )={\sf R}(\lambda -\frac{1}{2}\eta )$,
where ${\sf R}$ is the 8-vertex model $R$-matrix (\ref{02}).

As is proved in \cite{z},
the space $\Theta_{4\ell}^{+}$ for
$\ell \in \frac{1}{2}\ZZ_{+}$ is annihilated by the operator
\beq
\label{W}
{\bf w}_{\ell} = \sum_{k=0}^{2\ell +1} (-1)^k
\left [ \begin{array}{c} 2\ell +1 \\ k \end{array} \right ]
\, \frac{ \theta _{1}(2z+
2(2\ell -2k +1)\eta )}{\prod_{j=0}^{2\ell +1}
\theta_1(2z+2(j-k)\eta )} \,
e^{(2\ell -2k +1)\eta \p_z }.
\eeq
Here and below we use the
"elliptic factorial" and "elliptic binomial" notation:
\beq
\label{binom} [j]\equiv \theta_1(2j\eta)\,,
\;\;\;\;\;\; [n]!=\prod_{j=1}^{n}[j]\,,
\;\;\;\;\;\;
\left [ \begin{array}{c}n\\m\end{array}\right ]
\equiv \displaystyle{\frac{[n]!}{[m]![n-m]!}}\,.
\eeq
Another important property of the operator ${\bf w}_{\ell}$
established in \cite{z} is that ${\bf w}_{\ell}$
intertwines representations
of spin $\ell$ and of spin $-(\ell +1)$:
\beq
\label{S3}
{\bf s}_{a}^{(-\ell -1)} {\bf w}_{\ell} = {\bf w}_{\ell}
{\bf s}_{a}^{(\ell )}\,,
\;\;\;\;\;\; a=0,\ldots , 3\,.
\eeq
The same intertwining relation can be written
for the quantum $\mbox{{\sf L}}$-operator (\ref{L}).
The operator ${\bf w}_{\ell}$ is an elliptic analog of
$(d/dz)^{2\ell +1}$ in the following sense.
In the case of the algebra $sl_2$, the intertwining
operator between representations of spins $\ell$ and
$-\ell -1$ (realized by differential operators in $z$) is
$(d/dz)^{2\ell +1}$. It annihilates the linear space
of polynomials of degree $\leq 2\ell$.

The operator ${\bf w}_{\ell}$ can be extended to
arbitrary complex values of $\ell$
in which case it is represented by a half-infinite series in the
shift operator $e^{2\eta \p_z}$.
The series is an elliptic analog of the very-well-poised
basic hypergeometric series with an operator argument.
For its explicit form see (\ref{winf}).
The intertwining relations
(\ref{S3}) hold true in this more general case, too.
The proof is outlined in Appendix B.

In what follows we need
infinite dimensional representations of the Sklyanin
algebra in some functional space. Operators (\ref{Sa})
provide such a representation. However, they are not
sufficient for our purpose since it is necessary to
specify the functional space where these operators
are going to act.
It is not clear how to do that for continuous functions.
On the other hand, the difference character of the operators
(\ref{Sa}) allows one to restrict them to
{\it discontinuous} functions of the form described below.

Let $\delta (z)$ be the function equal to zero everywhere
but at $z=0$, where it equals $1$:
$\delta (z)=0$, $z\neq 0$, $\delta (0)=1$.
Clearly, $z\delta (z)=0$.
The Sklyanin algebra realized as in (\ref{Sa}) acts,
for any $\nu \in \CC$, in
the space of functions of the form
\beq
\label{inf1}
f(z)=\sum_{k\in\z}f_k \delta (z-\nu +2k\eta )\,,
\;\;\;\;\;\;\; f_k \in \CC \,.
\eeq
This space is isomorphic to the linear space of
sequences $\{ f_k \}_{k\in {\bf \z}}$. We call functions of the
form (\ref{inf1}) {\it combs}.

A comb is said to be finite from the right
(respectively, from the left)
if there exists $M\in \ZZ$ such that $f_k =0$ as $k>M$
(respectively, $k<M$). Let ${\cal C}^{\vdash}$ (respectively,
${\cal C}^{\dashv}$) be the space of combs finite from the left
(respectively, from the right).

We define the pairing
\beq
\label{inf2}
(F(z),\, \delta (z-a))=F(a)
\eeq
for any function $F(z)$, not necessarily of the form (\ref{inf1}).
In particular,
\beq
\label{inf3}
(\delta (z-a),\, \delta (z-b))=
\delta (a-b)\,.
\eeq
Formally, this pairing can be written as an integral:
$$
(F(z),\, \delta (z-a))=
\int F(z)\delta (z-a)dz
$$
(perhaps a $q$-integral symbol would be more appropriate).
We stress that the integral here means nothing more
than another notation for the pairing,
especially convenient in case of many variables.
By linearity, the pairing can be extended to the whole
space of combs. We note that the pairing between the spaces
${\cal C}^{\vdash}$ and ${\cal C}^{\dashv}$ is well defined
since the sum is always finite.

Combs are to be viewed as kernels of difference operators.
By a difference operator in one variable we mean any
expression of the form
\beq
\label{D1}
{\bf D}=\sum_{k\in {\z}}c_k(z)e^{(\mu +2k\eta )\p _{z}}\,,
\;\;\;\;\;\; \mu \in \CC\,.
\eeq
The comb
\beq
\label{D2}
D(z, \zeta )=
\sum_{k\in {\z}}c_k(z)
\delta (z-\zeta +\mu +2k\eta )\,,
\eeq
regarded as a function of any one of the variables $z$, $\zeta$,
is the kernel of this difference operator in the following sense.
Using the pairing introduced above, we can write:
\beq
\label{D3}
({\bf D}f)(z)=
\int D(z,\zeta )f(\zeta )d\zeta =
\sum_{k\in {\z}}c_k(z)f(z+\mu +2k\eta )\,.
\eeq
The kernel $D(z, \zeta )$ can be viewed as an infinite matrix
with continuously numbered rows ($z$) and columns ($\zeta$).
Then the convolution with respect to the second
argument of the kernel, as in (\ref{D3}),
defines action of the operator from the left.
The convolution with respect to the first argument
defines the action from the right,
\beq
\label{D4}
(f{\bf D})(z)=
\int f(\zeta )D(\zeta ,z )d\zeta \,,
\eeq
equivalent to the action of the transposed
difference operator from the left:
\beq
\label{D5}
{\bf D}^{{\sf t}}=
\sum_{k\in {\z}}
e^{-(\mu +2k\eta )\p _{z}}
c_k(z) =
\sum_{k\in {\z}}
c_k(z-\mu -2k\eta )
e^{-(\mu +2k\eta )\p _{z}}\,.
\eeq
The transposition ${\sf t}$ is the anti-automorphism of the
algebra of difference operators such that
$\bigl ( c(z)e^{\alpha \p_{z}}\bigr )^{{\sf t}}=
e^{-\alpha \p_{z}}c(z)$.
In terms of the above pairing we can write
$(f, {\bf D}g)=({\bf D}^{{\sf t}}f, g)$.

The kernels of Sklyanin's operators (\ref{Sa}) are:
\beq
\label{D7}
s_a(z,z')=
\frac{\theta _{a+1}(2z -2\ell \eta)}
{\theta _{1}(2z)}\,\delta (z-z' +\eta)
-\frac{\theta _{a+1}(-2z -2\ell \eta)}
{\theta _{1}(2z)}\, \delta (z-z' -\eta )\,.
\eeq
Let us also note that
${\bf s}_{a}^{(\ell){\sf t}}=-(-1)^{\delta_{a,0}}
\theta_{1}(2z)
{\bf s}_{a}^{(-\ell -1)}
(\theta_{1}(2z))^{-1}$.

\section{Local vacuum vectors}

We are going to adopt the Baxter method of vacuum vectors
to the infinite-dimensional representations
of the Sklyanin algebra described in the previous section.

Let us recall the general definition of vacuum vectors.
Consider an arbitrary ${\sf L}$-operator ${\sf L}$
with two-dimensional auxiliary space $\CC^2$, i.e., an arbitrary
$2\times 2$
operator-valued matrix whose matrix
elements are operators
${\bf L}_{11}$, ${\bf L}_{12}$, ${\bf L}_{21}$,
${\bf L}_{22}$:
\beq
{\sf L}=\left ( \begin{array}{cc} {\bf L}_{11}& {\bf L}_{12}\\
{\bf L}_{21}& {\bf L}_{22} \end{array} \right ).
\label{5.0}
\eeq
They act in
a linear space ${\cal H}$
which is called the {\it quantum space} of
the ${\sf L}$-operator.

For a moment, let $\psi$, $\phi$ denote vectors from
$\CC^2$.
Let $X\in {\cal H}$, then acting by the quantum
${\sf L}$-operator on the tensor product
$X\otimes \phi$, we, generally speaking, obtain a mixed state
in the quantum space:
${\sf L}X\otimes \phi =X_1\otimes \phi_1+X_2\otimes \phi_2$.
The special case of a pure state,
\beq
\label{V1}
{\sf L}X\otimes \phi =X'\otimes \psi\,,
\eeq
is of particular importance. Given an ${\sf L}$-operator,
Baxter's receipe is to find all pairs $X, \phi$ such that
(\ref{V1}) holds. Taking the scalar product
with the vector $\psi^{\bot}$ orthogonal to $\psi$, we get:
\beq
\label{V2}
\bigl (\psi^{\bot}{\sf L}\phi \bigr )X=0\,,
\eeq
i.e., the operator
${\bf K}=\bigl (\psi^{\bot}{\sf L}\phi \bigr )$ (acting in the
quantum space only) has a zero mode
$X\in {\cal H}$.
Suppose (\ref{V1}) (or (\ref{V2})) holds with some
vectors $\phi$, $\psi$;
then the vector $X$ is called
a {\it vacuum vector} of the $L$-operator.

A few remarks are in order.
The relation (\ref{V1}) (in the particular case ${\cal H}\cong
\CC^2$) was the starting point for Baxter in his solution of
the 8-vertex and $XYZ$ models \cite{Baxter}. This is
what he called the ``pair-propagation through a vertex"
property. In the context of the
quantum inverse scattering method \cite{FT}
the equivalent condition (\ref{V2}) is more customary.
It defines local vacua of the
(gauge-transformed) ${\sf L}$-operator.
A generalization of that
solution to the higher spin $XYZ$ model was given
by Takebe in \cite{Takebe}. An algebro-geometric
approach to the equation (\ref{V1}) for finite-dimensional
matrices ${\bf L}_{ik}$ was developed by Krichever \cite{krivac}.

For ${\sf L}$-operators with elliptic spectral parameter
it is convenient to pass to the elliptic parametrization
of the components of the vectors $\phi$, $\psi$.
We introduce the 2-component vector
$$
\bigl |\zeta \bigr >=
\left (
\begin{array}{l}\bar \theta_{4}(\zeta)\\
\bar \theta_{3}(\zeta)
\end{array}
\right )
$$
and the corresponding covector
$\bigl < \zeta \bigr |= \bigl (
\bar \theta_{4}(\zeta),\,
\bar \theta_{3}(\zeta) \bigr )$
(recall that
$\bar \theta_{a}(\zeta)= \theta _{a}(\zeta |\frac{\tau}{2})$).
The vector orthogonal to $\bigl < \zeta \bigr |$ is
$\bigl |\zeta \bigr >^{\bot}=
\left (
\begin{array}{r}\bar \theta_{3}(\zeta)\\
-\bar \theta_{4}(\zeta)
\end{array}
\right )$,
$\bigl < \zeta \bigr | \zeta \bigr >^{\bot}=0$.
More generally, we have:
\beq
\label{scprod}
\bigl < \xi \bigr |\zeta \bigr >^{\bot}=
2\theta_{1}(\xi +\zeta )\theta_{1}(\xi -\zeta )\,.
\eeq
Note also that
\beq
\label{orth}
\bigl | \zeta +{\scriptstyle \frac{1}{2} }
(1+\tau )\bigr >=
e^{-\frac{\pi i \tau}{2} -2\pi i \zeta}
\bigl | \zeta \bigr >^{\bot}.
\eeq

Let us specify the general definition of the vacuum
vectors for the elliptic ${\sf L}$-operator
(\ref{L}). Writing ${\sf L}(\lambda )\bigl |\zeta \bigr >$
(respectively, $\bigl <\zeta \bigr |{\sf L}(\lambda)$)
we mean that the 2$\times$2 matrix (\ref{L}) acts on the
2-component vector from the left (respectively, on the
2-component covector from the right), as usual.
Similarly, we introduce right and left vacuum vectors
$X_R$, $X_L$ according to the relations
\beq
\label{V3}
\bigl < \zeta \bigr | {\sf L}(\lambda)X_R =
\bigl < \xi \bigr | X'_R \,,
\;\;\;\;\;\;
X_L\bigl < \zeta \bigr | {\sf L}(\lambda) =
X'_L \bigl < \xi \bigr |\,.
\eeq
In the latter formula the difference
operators ${\bf s}_a$ act on $X_L$ from the right.

The following pictorial representation is helpful:

\begin{center}
\special{em:linewidth 0.4pt}
\unitlength 0.4mm
\linethickness{0.4pt}
\begin{picture}(91.00,55.00)
\emline{29.67}{30.00}{1}{80.33}{30.00}{2}
\special{em:linewidth 1.6pt}
\emline{55.00}{5.00}{3}{55.00}{55.00}{4}
\put(55.00,0.00){\makebox(0,0)[cc]{$X(z)$}}
\put(19.67,30.00){\makebox(0,0)[cc]{$\bigl <\zeta \bigr |$}}
\put(91.00,30.00){\makebox(0,0)[cc]{$\bigl |\xi \bigr >^{\bot}$}}
\end{picture}
\end{center}

\noindent
Here we use the standard conventions of the quantum inverse
scattering method.

Introducing the operator
\beq
\label{V5}
{\bf K}= {\bf K}(\zeta , \xi )=
\bigl < \zeta \bigr | {\sf L}(\lambda)\bigl |\xi \bigr >^{\bot},
\eeq
we can rewrite (\ref{V3}) as
${\bf K}X_R =X_L{\bf K}=0$.
The explicit form of the operator ${\bf K}$ is easily
found from (\ref{fact}):
\beq
\label{E1}
{\bf K}={\bf K}(\zeta, \xi)=\rho (z) e^{\eta \p_{z}}+
\rho (-z) e^{-\eta \p_{z}}\,,
\eeq
where
$$
\rho (z)=\frac{1}{\theta_{1}(2z)}
\prod_{\epsilon =\pm}\theta_{1}(z+\epsilon \zeta
-\lambda_{+} )
\theta_{1}(z+\epsilon \xi
+\lambda_{-} )\,.
$$
So, the equations for the right and left vacuum vectors
read
\beq
\label{R1}
\rho (z)X_{R}(z+\eta)=
-\rho (-z)X_{R}(z-\eta)\,,
\eeq
\beq
\label{L1}
\rho (-z-\eta)X_{L}(z+\eta)=
-\rho (z-\eta)X_{L}(z-\eta)\,.
\eeq

These equations themselves are not yet enough.
An important additional input is the space where to solve
them. Our choice, justified by the final results,
is the following. We are going to solve
the first equation in the space
of combs finite from the left,
${\cal C}^{\vdash}$, while the second equation
is to be solved in the space of combs finite from the
right, ${\cal C}^{\dashv}$. (For precise definitions
of these spaces see Sec.\,3.)

Let us begin by solving equations (\ref{R1}), (\ref{L1}) in the
space of meromorphic functions of $z$. This is easy to
do in terms of the elliptic gamma-function
$\Gamma (z|\tau, 2\eta)$ introduced by Ruijsenaars \cite{R3}
and recently studied in detail by
Felder and Varchenko \cite{FV3}. This function was
introduced just to solve equations of such a type.
Its main properties are summarized in Appendix A.
In what follows we write
$\Gamma (z|\tau, 2\eta)=\Gamma(z)$ for brevity.
Let ${\cal X}_{R}(z)$,
${\cal X}_{L}(z)$ be meromorphic solutions
of eqs.\,(\ref{R1}), (\ref{L1}), respectively.
Using (\ref{gamma3}), we easily find:
\beq
\label{R2}
{\cal X}_{R}(z)=e^{-4\pi i \ell z}
\omega (\lambda , \zeta \!-\!\xi )
\prod_{\epsilon =\pm}
\frac{\Gamma (z+\epsilon \zeta +\lambda_{+} +\eta )
\Gamma (z+\epsilon \xi -\lambda_{-} +\eta )}
{\Gamma (z+\epsilon \zeta -\lambda_{+} +\eta )
\Gamma (z+\epsilon \xi +\lambda_{-} +\eta )}\,,
\eeq
\beq
\label{L2}
{\cal X}_{L}(z)=
e^{4\pi i (\ell +1)z}
\Phi_{\ell}(2\lambda)
\theta_{1}(2z)
\prod_{\epsilon =\pm}
\frac{\Gamma (z+\epsilon \zeta -\lambda_{+} )
\Gamma (z+\epsilon \xi +\lambda_{-} )}
{\Gamma (z+\epsilon \zeta +\lambda_{+} +2\eta )
\Gamma (z+\epsilon \xi -\lambda_{-} +2\eta )}\,.
\eeq
The $z$-independent normalization factors
$\omega (\lambda , \zeta \!-\!\xi )$ in (\ref{R2})
and $\Phi_{\ell}(2\lambda)$
in (\ref{L2}) are introduced here for later convenience.
The function
$\omega (\lambda , \zeta )$ is given by
\beq
\label{omega}
\omega (\lambda , \zeta )=
\,
\frac{\theta_1(2\lambda \!-\!2\ell \eta |\,2\eta )}
{\theta_1 (2\lambda +\zeta |\,2\eta \,)}
\, \exp \left ( -\, \frac{\pi i \zeta^2}{2\eta }
-4\pi i \ell \lambda \right )
\eeq
(the $\theta$-functions here have the modular parameter
$2\eta$ rather than $\tau$!).
The function $\Phi_{\ell}$
obeys the difference equation
$\theta_1 (z-2\ell \eta )\Phi_{\ell}(z+2\eta )=
\theta_1 (z+2(\ell \!+\!1) \eta )\Phi_{\ell}(z)$.
For $\ell \in \frac{1}{2}\ZZ_{+}$ the solution is very simple:
\beq
\label{Phi1}
\Phi_{\ell}(z)=
\prod_{k=-\ell}^{\ell}\theta_{1}(z +2k\eta)\,,
\eeq
where $k$ ranges over
$-\ell, \, -\ell +1 , \, \ldots , \ell$.
For arbitrary complex values of $\ell$
the function $\Phi_{\ell}$ can be defined
through the elliptic gamma-function (see (\ref{Phig})).
If necessary, we write
${\cal X}_{R,L}(z)={\cal X}_{R,L}(z;\zeta , \xi , \lambda)$.
Note that
${\cal X}_{R,L}(z;\zeta , \xi , \lambda)$
differs from
${\cal X}_{R,L}(z;\xi , \zeta , -\lambda)$
only by a $z$-independent factor.

An important remark is in order. The solutions
(\ref{R2}), (\ref{L2}) are by no means unique.
Clearly, one can multiply them by an arbitrary function
$f(z)$ such that $f(z+\eta )=f(z-\eta )$. In this paper,
we do not take care of such $2\eta$-periodic factors
fixing them to be constant. In other words,
(\ref{R2}), (\ref{L2}) are ``minimal" solutions
to eqs.\,(\ref{R1}), (\ref{L1}) in the sense of \cite{R3}.
By this we break the symmetry
between the three periods $1, 2\eta , \tau$ involved in the
problem. In particular, the modular symmetry is broken since
$2\eta$-periodic factors may be important
for good modular properties of the solutions. Indeed,
the modular transformation of the elliptic gamma-function
has the form \cite{FV3}:
$$
\Gamma (z|\tau , 2\eta )=e^{i\pi P(z)}
\frac{\Gamma (\frac{z}{\tau}|-\frac{1}{\tau}, \,\frac{2\eta}{\tau})}
{\Gamma (\frac{z-\tau}{2\eta}|-\frac{1}{2\eta}, \,
\frac{\tau}{2\eta})}\,,
$$
where $P(z)$ is a polynomial of third degree (see (\ref{polP})).
The denominator
is a $2\eta$-periodic function which does not enter the
minimal solution of the modular transformed equation.
We note in passing that this formula is
an elliptic version of eq.\,(40) from the paper \cite{modular},
where a modular double of quantum group was introduced.
We hope to discuss this point elsewhere. See also the end
of Sec.\,8.

Now we can proceed to built a half-infinite comb solution.
Any infinite comb of the form
$X(z)={\cal X}_{R,L}(z)\sum_{k\in {\z}}\delta (z-z_0 +2k\eta)$
automatically satisfies eq.\,(\ref{R1}) (resp., (\ref{L1})).
How to truncate the
comb holding these equations valid?
Let ${\cal X}(z)$ be any one of ${\cal X}_R(z)$, ${\cal X}_L(z)$.
The position of zeros
and poles of the $\Gamma$-function suggests two possibilities:
\begin{itemize}
\item[(1)]
To set $z_0$ to be a zero of the function
${\cal X}$: ${\cal X}(z_0)=0$. If this zero comes from
a pole of one of the $\Gamma$-functions downstairs, then
${\cal X}(z_0-2k\eta)=0$, $k\geq 0$, and the comb is truncating
from the left. If $z_0$ is a zero of
one of the $\Gamma$-functions
upstairs, then
${\cal X}(z_0-2k\eta)=0$, $k<0$, and the comb is truncating
from the right.
\item[(2)]
To set $z_0$ to be a pole of ${\cal X}$ and define
$X(z)=\sum_{k\in \z} \,\mbox{res}\Bigl |_{z=z_0-2k\eta}\!\!
{\cal X}(z)\delta (z-z_0 +2k\eta)$. If $z_0$ is a pole
of a $\Gamma$-function upstairs, then the comb is truncating
from the right. If the pole comes from a zero of one of the
$\Gamma$-functions downstairs, then the comb is truncating
from the left.
\end{itemize}
In the sequel we use both possibilities:
(1) for $X_R$ and (2) for $X_L$.
Other options are to be further analysed.

Assume that the parameters are in general position,
i.e., there is no cancellation of zeros and poles
in (\ref{R2}). (For $\ell \in \frac{1}{2}\ZZ _{+}$ the
sufficient conditions are $\lambda \neq 2m\eta$,
$\xi -\zeta \neq 2m\eta$, $m\in \ZZ$.)
Then ${\cal X}_{R}(\xi -\lambda_{-} -\eta -2k\eta )=0$
for all integer $k\geq 0$. At the same time this function
is not identically zero when $k$ is negative. Therefore,
\beq
\label{R3}
%X_R^+(z)=
X_R^+(z;\zeta , \xi , \lambda)=
{\cal X}_{R}(z; \zeta , \xi , \lambda)
\sum_{k=1}^{\infty}
\delta (z-\xi +\lambda_{-} +\eta -2k\eta )
\eeq
is a non-trivial comb finite from the left
which solves eq.\,(\ref{R1}). Clearly, this is not the only
possibility to truncate the infinite comb from the left,
even if to stick to (1). There are three
other series of zeros of the function (\ref{R2}) coming from
poles of the gamma-functions in the denominator.
Any one of them can be used to generate, in the same way,
a half-infinite comb solution that belongs to ${\cal C}^{\vdash}$.
However, inasmuch as the vacuum vectors are going to be
kernels of difference operators,
we want the $\delta$-functions
in the vacuum vectors to depend on the differences
$z-\zeta$ or $z-\xi$. So, along with (\ref{R3}), we have
another comb obeying eq.\,(\ref{R1}):
\beq
\label{R3a}
%X_R^-(z)=
X_R^-(z; \zeta ,\xi , \lambda)=
{\cal X}_{R}(z; \xi , \zeta, -\lambda)
\sum_{k=1}^{\infty}
\delta (z-\zeta -\lambda_{+} +\eta -2k\eta )\,.
\eeq

For $X_L$ we use (2):
${\cal X}_{L}(z)$ has poles at
$z=\zeta +\lambda_{+} -2k\eta$, $k\geq 0$.
In general position, these poles are all simple.
The lattice of poles is finite from the right.
Therefore, we can construct a solution to eq.\,(\ref{L1})
in the space of combs finite from the right as follows:
\beq
\label{L3}
%X_L^+(z)=
X_L^+(\lambda, \zeta , \xi ;z)=
\sum_{k=0}^{\infty}
\,\mbox{res}\,\Bigl |_{z=\zeta +\lambda_{+} -2k\eta}
\!\!\!{\cal X}_{L}(z;\zeta , \xi , \lambda)\,\,
\delta (z-\zeta -\lambda_{+} +2k\eta)
\eeq
Similarly to the case of the right vacuum vectors,
we have a second solution:
\beq
\label{L3a}
%X_L^-(z)=
X_L^-(\lambda, \zeta , \xi ;z)=
\sum_{k=0}^{\infty}
\,\mbox{res}\,\Bigl |_{z=\xi -\lambda_{-} -2k\eta}
\!\!\!{\cal X}_{L}(z;\xi , \zeta , -\lambda)\,\,
\delta (z-\xi +\lambda_{-} +2k\eta)\,.
\eeq

To make these formulas closer to (\ref{R3}), (\ref{R3a}),
we introduce the auxiliary function
$$
{\cal X}_{L}^{\varepsilon}(z;\zeta , \xi , \lambda )=
e^{4\pi i (\ell +1)z}
\Phi_{\ell}(2\lambda)
\theta_{1}(2z)  \times
$$
$$
\times \frac{\Gamma (z -\zeta -\lambda_{+}+\varepsilon )
\Gamma (z +\zeta -\lambda_{+} )
\Gamma (z - \xi +\lambda_{-}+\varepsilon )
\Gamma (z + \xi +\lambda_{-} )}{\Gamma
(z -\zeta +\lambda_{+}+2\eta )
\Gamma (z +\zeta +\lambda_{+}+2\eta )
\Gamma (z - \xi -\lambda_{-}+2\eta )
\Gamma (z + \xi -\lambda_{-} +2\eta )}\,,
$$
then, up to an irrelevant constant factor,
\beq
\label{L31}
X_L^+(z;\zeta ,\xi ,\lambda )=
\lim_{\varepsilon \to 0}\,
\varepsilon
{\cal X}_{L}^{\varepsilon}(z;\zeta , \xi , \lambda )
\sum_{k=0}^{\infty}
\delta (z-\zeta -\lambda_{+} +2k\eta)\,,
\eeq
\beq
\label{L3a1}
X_L^-(\lambda , \zeta , \xi ;z)=
\lim_{\varepsilon \to 0}\,
\varepsilon
{\cal X}_{L}^{\varepsilon}(z; \xi , \zeta , -\lambda )
\sum_{k=0}^{\infty}
\delta (z-\xi +\lambda_{-} +2k\eta)\,.
\eeq
Here one should first substitute the values of $z$
in the $\delta$-functions and then take the limit.

Note that the function $\rho (z)$ and the meromorphic
solutions ${\cal X}_{R,L}$ are invariant with respect to
the simultaneous change $\lambda
\rightarrow -\lambda$, $\zeta \leftrightarrow \xi$ while the
half-infinite combs $X_{R,L}^{\pm}$ are
transformed one into another:
\beq
\label{E2}
\begin{array}{l}
X_R^-(z;\zeta , \xi , \lambda)=
%\displaystyle{\frac{\omega (\lambda , \zeta -\xi )}
%{\omega (-\lambda , \xi -\zeta )}}
X_R^+(z;\xi , \zeta , -\lambda)\,,\\ \\
X_L^-(\lambda , \zeta , \xi ;z)=
%\displaystyle{\frac{\Phi_{\ell}(2\lambda )}
%{\Phi_{\ell}(-2\lambda )}}
X_L^+(-\lambda , \xi , \zeta ;z)\,.
\end{array}
\eeq

The construction of local vacuum vectors is
completed. The vacuum vectors are defined up to an
arbitrary normalization factor which may depend on
$\zeta, \xi$ and $\lambda$.
In our normalization, the final formulas
have the most simple form.

By a straightforward
computation, which uses some standard identities
for theta-functions,
one can prove the following relations valid for
$X_{R,L}^{\pm}=X_{R,L}$:
\beq
\label{R4}
\bigl < \zeta \bigr | {\sf L}(\lambda )
X_{R}(z;\zeta , \xi , \lambda )
=\theta_{1}(2\lambda_{-})
\bigl < \xi \bigr |
X_{R}(z;\zeta , \xi , \lambda +\eta )\,,
\eeq
\beq
\label{L4}
X_{L}(\zeta , \xi , \lambda ;z )
\bigl < \zeta \bigr | {\sf L}(\lambda )
=\theta_{1}(2\lambda_{+})
X_{L}(\zeta , \xi , \lambda -\eta ;z)
\bigl < \xi \bigr | \,,
\eeq
which are simply (\ref{V3}) in the specified parametrization.
In the first equation, each operator entry of the matrix
${\sf L}(\lambda )$ acts on the vacuum vector from the left
while in the second one they act from the right.
In the both equations, the
two-component covector $\bigl < \zeta \bigr |$ is applied
from the left.
Using
(\ref{trans}) and (\ref{orth}), one can represent
these equations in the form, where the two-component vector
is applied from the right:
\beq
\label{R5}
{\sf L}(\lambda )\bigl | \xi \bigr >^{\bot}
X_{R}(z;\zeta , \xi , \lambda )
=\theta_{1}(2\lambda_{+} )
\bigl | \zeta \bigr >^{\bot}
X_{R}(z;\zeta , \xi , \lambda -\eta )\,,
\eeq
\beq
\label{L5}
X_{L}(\zeta , \xi , \lambda ;z )
{\sf L}(\lambda )\bigl | \xi \bigr >^{\bot}
=\theta_{1}(2\lambda_{-})
\bigl | \zeta \bigr >^{\bot}
X_{L}(\zeta , \xi , \lambda +\eta ;z )\,.
\eeq
For the holomorphic vacuum vectors
((\ref{R2}) at $\ell \in \ZZ _{+})$ these relations
hold \cite{Takebe} under the additional restriction
$\xi -\zeta \in 2\eta \ZZ$.

\section{Global vacuum vectors}

The next step, in accordance
with \cite{Baxter}, is to find the vacuum vectors for each
$\stackrel{i}{{\sf L}}(\lambda)$ and consider the action
of the transfer matrix ${\bf T}(\lambda)$ (\ref{T3}) on their
direct product along the chain. It then appears that the
result of this action has a very special and simple form:
it is a sum of two terms, each of them being a vacuum vector
that belongs to the same family. From the technical point of
view, the easiest way to see this is to use the gauge
transformation argument, now standard.

Let us make the ``gauge transformation" of the ${\sf L}$-operator:
\beq
\label{G1}
\stackrel{i}{{\sf L}}(\lambda)\, \longrightarrow  \,
G(\zeta _i)
\stackrel{i}{{\sf L}}(\lambda)
(G(\zeta _{i+1}))^{-1}=
\stackrel{\!i}{{\sf L}'}(\lambda)\,,
\eeq
where $G(\zeta)$ are some 2$\times$2 $c$-number matrices
depending on the parameter $\zeta$. If $\zeta_{N+1}=\zeta_1$,
then the trace is not changed:
$$
{\bf T}(\lambda)=
\mbox{tr}\, \left (
\stackrel{\!1}{{\sf L}'}(\lambda)
\stackrel{\!2}{{\sf L}'}(\lambda)
\,\ldots \,
\stackrel{\!N}{{\sf L}'}(\lambda)\right )\,.
$$
We set the matrix $G(\zeta )$ to be
\beq
\label{G2}
G(\zeta )=\left (
\begin{array}{cc}
\bar \theta_4 (\zeta )&
\bar \theta_3 (\zeta ) \\
0&1 \end{array} \right )\,.
\eeq
Then, acting on a vacuum vector $X_R$ by each element
of the gauge-transformed ${\sf L}$-operator and using
the relations (\ref{R4}), (\ref{R5}), we obtain:
$$
\left ( \begin{array}{cc} {\bf L}'_{11}X_R
& {\bf L}'_{12}X_R\\
{\bf L}'_{21}X_R& {\bf L}'_{22}X_R \end{array} \right )=
$$
$$
=\frac{1}{\bar \theta_4 (\xi)}
\left (
\begin{array}{cc}
\bar \theta_4 (\zeta )&
\bar \theta_3 (\zeta ) \\
0&1 \end{array} \right )
\left ( \begin{array}{cc} {\bf L}_{11}X_R
& {\bf L}_{12}X_R\\
{\bf L}_{21}X_R& {\bf L}_{22}X_R \end{array} \right )
\left (
\begin{array}{rr}1&
-\bar \theta_3 (\xi )\\
0&\bar \theta_4 (\xi )
\end{array} \right ) =
$$
\beq
\label{triang}
=\left (
\begin{array}{ccc}
X_R'&&0\\ &&\\
\ast  && \displaystyle{
\frac{\bar \theta_4 (\zeta)}{\bar \theta_4 (\xi)}}X_{R}''
\end{array} \right ).
\eeq
Here $X_{R}'=\theta_{1}(2\lambda_{-})
X_{R}(z_i;\zeta , \xi , \lambda +\eta )$,
$X_{R}''=\theta_{1}(2\lambda_{+})
X_{R}(z_i;\zeta , \xi , \lambda -\eta )$ and we imply
that $\zeta_i =\zeta$, $\zeta_{i+1}=\xi$.
The star in the left-down corner stands for something
irrelevant for us.

We define the global vacuum vectors
$\vec X_{R,L}^{\pm}$
to be products of the local ones over all sites of the chain:
\beq
\label{G3R}
\vec X^{\pm}_{R}(\lambda)=
\vec X^{\pm}_{R}(\vec z, \vec \zeta , \lambda)=
\prod_{i=1}^{N}
X_{R}^{\pm}(z_i; \zeta_i ,\zeta_{i+1}, \lambda )\,,
\eeq
\beq
\label{G3L}
\vec X^{\pm}_{L}(\lambda)=
\vec X^{\pm}_{L}(\lambda, \vec \zeta ,\vec z)=
\prod_{i=1}^{N}
X_{L}^{\pm}(\lambda ,\zeta_i ,\zeta_{i+1};z_i )\,,
\eeq
where the notation $\vec z =\{z_1 , \ldots , z_N\}$,
$\vec \zeta =\{\zeta_1 , \ldots , \zeta_N\}$ is used.
The parameters $\zeta_i$ are naturally asociated with
the horizontal edges of the lattice:

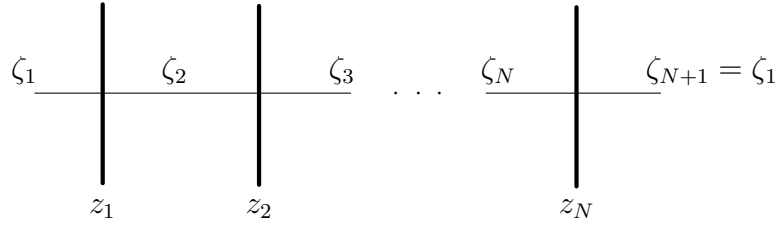
\begin{center}
\special{em:linewidth 0.4pt}
\unitlength 0.6mm
\linethickness{0.4pt}
\begin{picture}(158.67,49.67)
\emline{20.00}{30.00}{1}{90.00}{30.00}{2}
\put(100.00,30.00){\circle*{0.67}}
\put(105.00,30.00){\circle*{0.67}}
\put(109.67,30.00){\circle*{0.67}}
\emline{120.00}{30.00}{3}{158.67}{30.00}{4}
\special{em:linewidth 1.6pt}
\emline{35.00}{10.00}{5}{35.00}{49.67}{6}
\emline{69.67}{9.67}{7}{69.67}{49.33}{8}
\emline{140.00}{9.33}{9}{140.00}{49.00}{10}
\put(17.67,35.00){\makebox(0,0)[cc]{$\zeta _1$}}
\put(51.00,35.00){\makebox(0,0)[cc]{$\zeta _2$}}
\put(88.00,35.00){\makebox(0,0)[cc]{$\zeta _3$}}
\put(122.67,35.00){\makebox(0,0)[cc]{$\zeta _N$}}
\put(170.00,35.00){\makebox(0,0)[cc]{$\zeta _{N+1} = \zeta _1$}}
\put(35.00,4.33){\makebox(0,0)[cc]{$z_1$}}
\put(69.67,4.33){\makebox(0,0)[cc]{$z_2$}}
\put(140.00,4.33){\makebox(0,0)[cc]{$z_N$}}
\end{picture}
\end{center}

\noindent
Then
$$
{\bf T}(\lambda) \vec X^{\pm}_{R}(\lambda)=
\mbox{tr}\,\prod_{i=1}^{N}
\left (
\stackrel{\!i}{{\sf L}'}(\lambda)
X_{R}^{\pm}(z_i; \zeta_i ,\zeta_{i+1}, \lambda )\right )
$$
and similarly for the $\vec X^{\pm}_{L}(\lambda)$.
The triangularity of the matrix (\ref{triang})
allows one to compute the trace. Taking into account
the periodic boundary condition, we arrive at
the following relations:
\beq
\label{G4R}
{\bf T}(\lambda) \vec X^{\pm}_{R}(\lambda)=
\theta_{1}^{N}(2\lambda_{-} )
\vec X^{\pm}_{R}(\lambda +\eta)+
\theta_{1}^{N}(2\lambda_{+})
\vec X^{\pm}_{R}(\lambda -\eta)\,,
\eeq
\beq
\label{G4L}
\vec X^{\pm}_{L}(\lambda){\bf T}(\lambda )=
\theta_{1}^{N}(2\lambda_{-})
\vec X^{\pm}_{L}(\lambda +\eta)+
\theta_{1}^{N}(2\lambda_{+})
\vec X^{\pm}_{L}(\lambda -\eta)\,.
\eeq

Let us introduce left and right $Q$-operators,
${\bf Q}_{L}^{\pm}$ and
${\bf Q}_{R}^{\pm}$
(the {\it pre-$Q$-operators}).
They are difference operators
acting on functions $f(z_1 , \ldots , z_N)$ of $N$ variables,
respectively from the left and from the right, as follows:
\beq
\label{kQL}
({\bf Q}_{L}^{\pm}f)(\vec \zeta )=
\int d\vec z
\, \vec X^{\pm}_{L}(\lambda , \vec \zeta , \vec z)f(\vec z)\,,
\eeq
\beq
\label{kQR}
(f{\bf Q}_{R}^{\pm})(\vec \zeta )=
\int d\vec z
\, f(\vec z) \vec X^{\pm}_{R}(\vec z , \vec \zeta , \lambda )\,.
\eeq
The kernels here are
the global vacuum vectors defined above.
The symbol $\int d\vec z =\int \ldots \int dz_1 \ldots dz_N$ means
the natural pairing in the tensor product
of the local functional spaces.
On each factor, this pairing coincides with (\ref{inf2}).

The formulas (\ref{G4R}), (\ref{G4L}) can be written as
the operator relations:
\beq
\label{TQL}
{\bf T}(\lambda) {\bf Q}^{\pm}_{R}(\lambda)=
\theta_{1}^{N}(2\lambda_{-} )
{\bf Q}^{\pm}_{R}(\lambda +\eta)+
\theta_{1}^{N}(2\lambda_{+})
{\bf Q}^{\pm}_{R}(\lambda -\eta)\,,
\eeq
\beq
\label{TQR}
{\bf Q}^{\pm}_{L}(\lambda) {\bf T}(\lambda)=
\theta_{1}^{N}(2\lambda_{-} )
{\bf Q}^{\pm}_{L}(\lambda +\eta)+
\theta_{1}^{N}(2\lambda_{+})
{\bf Q}^{\pm}_{L}(\lambda -\eta)\,.
\eeq
One might multiply the pre-$Q$-operators by a scalar
normalization factor $f(\lambda )$ such that
$f(\lambda ) =f(\lambda +\eta )$. Clearly,
these relations remain unchanged.

\section{Basic commutation relations}

In this section we prove the basic commutation relations
\beq
\label{Q1}
{\bf Q}^{\pm}_L(\lambda ){\bf Q}^{\pm}_R(\lambda ')=
{\bf Q}^{\pm}_L(\lambda '){\bf Q}^{\pm}_R(\lambda )\,,
\eeq
\beq
\label{Q1a}
{\bf Q}^{+}_L(\lambda ){\bf Q}^{-}_R(\lambda ')=
{\bf Q}^{-}_L(\lambda '){\bf Q}^{+}_R(\lambda )\,.
\eeq

Let us start from the proof of (\ref{Q1}) for ${\bf Q}^{+}_{L,R}$.
To this end, consider the function
\beq
\label{Q2}
\int d^N \vec z \,\vec X_L^+ (\lambda ,\vec \zeta , \vec z)
\vec X_R^+ (\vec z , \vec \zeta ', \lambda ')
=\prod_{i=1}^{N}
\int dz_i X_L^+ (\lambda , \zeta_i , \zeta_{i+1};z_i )
X_R^+ (z_i ; \zeta'_i , \zeta'_{i+1}, \lambda ')
\eeq
which is the kernel of the operator product in the
left hand side of (\ref{Q1}).
Schematically, the structure of the kernel
can be illustrated as follows:

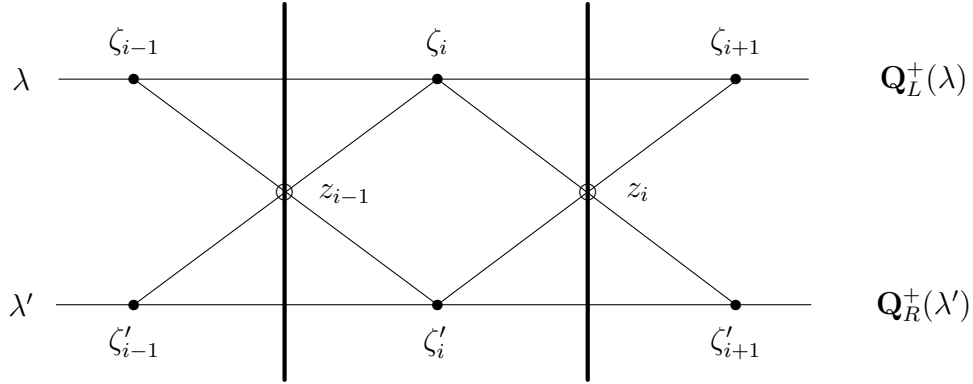
\begin{center}
\special{em:linewidth 0.4pt}
\unitlength 1mm
\linethickness{0.4pt}
\begin{picture}(125.00,50.00)
\emline{9.67}{10.00}{1}{110.00}{10.00}{2}
\emline{10.00}{40.00}{3}{109.67}{40.00}{4}
\put(60.33,40.00){\circle*{1.33}}
\put(60.33,10.00){\circle*{1.33}}
\put(100.00,10.00){\circle*{1.33}}
\put(100.00,40.00){\circle*{1.33}}
\put(20.00,40.00){\circle*{1.33}}
\put(20.00,10.00){\circle*{1.33}}
\emline{20.00}{10.00}{5}{60.33}{40.00}{6}
\emline{60.33}{40.00}{7}{100.00}{10.00}{8}
\emline{20.00}{40.00}{9}{60.33}{10.00}{10}
\emline{60.33}{10.00}{11}{100.33}{40.00}{12}
\put(40.00,25.00){\circle{2.00}}
\put(80.33,25.00){\circle{2.00}}
\put(5.00,40.00){\makebox(0,0)[cc]{$\lambda$}}
\put(5.00,10.00){\makebox(0,0)[cc]{$\lambda '$}}
\put(20.00,45.00){\makebox(0,0)[cc]{$\zeta _{i-1}$}}
\put(20.00,5.00){\makebox(0,0)[cc]{$\zeta '_{i-1}$}}
\put(60.33,45.00){\makebox(0,0)[cc]{$\zeta _i$}}
\put(60.00,5.00){\makebox(0,0)[cc]{$\zeta '_i$}}
\put(100.00,45.00){\makebox(0,0)[cc]{$\zeta _{i+1}$}}
\put(100.00,5.00){\makebox(0,0)[cc]{$\zeta '_{i+1}$}}
\special{em:linewidth 1.6pt}
\emline{40.00}{50.00}{13}{40.00}{0.00}{14}
\emline{80.33}{50.00}{15}{80.33}{0.00}{16}
\put(48.00,24.67){\makebox(0,0)[cc]{$z_{i-1}$}}
\put(87.00,25.00){\makebox(0,0)[cc]{$z_i$}}
\put(125.00,40.00){\makebox(0,0)[cc]{${\bf Q}^{+}_L(\lambda)$}}
\put(125.00,10.00){\makebox(0,0)[cc]{${\bf Q}^{+}_R(\lambda ')$}}
\end{picture}
\end{center}

\noindent
As is clear from (\ref{Q2}),
the kernel is a product of local multipliers associated
with each pair $(i, i\!+\!1)$ of adjacent sites of the chain.
For a while, we set $\zeta_i =\zeta$, $\zeta_{i+1}=\xi$.
Therefore, we should investigate the function
\beq
\label{Q3}
F^{\zeta \xi}_{\zeta ' \xi '}(\lambda , \lambda ')=
\int dz X_{L}^{+} (\lambda , \zeta , \xi ;z)
X_{R}^{+} (z; \zeta' , \xi ', \lambda ')\,.
\eeq

Let us begin with the $X_{L}^{+}$. Up to an irrelevant constant
factor, (\ref{L31}) reads
\beq
\label{Q4}
X_{L}^{+} (\lambda , \zeta , \xi ;z)=
\lim _{\varepsilon \to 0}\Bigl [ \varepsilon
{\cal X}^{\varepsilon}_{L}( \zeta +\lambda_{+}; \zeta , \xi ,
\lambda )\Bigr ]
\sum _{k\geq 0}x_{k}(\zeta , \xi , \lambda )
\delta ( z-\zeta -\lambda _{+} +2k\eta )\,,
\eeq
where
\beq
\label{xk}
\begin{array}{ll}
&x_{k}(\zeta , \xi , \lambda ) =
\displaystyle{ \frac{1}{[k]!} \,
\frac{\theta_1 (2\zeta \!+\!2\lambda_{+}\!-\!4k\eta )}
{\theta_1 (2\zeta \!+\!2\lambda_{+})} \, \times }
\\&\\
\times &\displaystyle{
\prod_{j=0}^{k-1}\frac{\theta_1(-2\lambda_{+}\!+\!2j\eta )
\theta_1(-2\lambda_{+}\!-\!2\zeta \!+\!2j\eta )
\theta_1(\xi \!-\!\zeta \!-\!2\ell \eta \!+\!2j\eta )
\theta_1(-\xi \!-\!\zeta \!-\!2\ell \eta \!+\!2j\eta )}
{\theta_1(-\!2\zeta \!+\!2(j\!+\!1)\eta )
\theta_1(\xi \!-\!\zeta \!-\!2\lambda \!+\!2(j\!+\!1)\eta )
\theta_1(-\xi \!-\!\zeta \!-\!2\lambda \!+\!2(j\!+\!1)\eta )} }\,.
\end{array}
\eeq
Substituting (\ref{Q4}) and
$
X_{R}^{+} (z; \zeta' , \xi ', \lambda ')=
{\cal X}_{R}(z; \zeta ' , \xi ' , \lambda ')
\sum_{k'\geq 0} \delta
(z\!-\!\xi ' \!+\!\lambda'_{-}\!-\!\eta \!-\!2k'\eta )
$
(see (\ref{R3})) into (\ref{Q3}),
performing the $z$-integral according to the
pairing rule (\ref{inf3}),
and passing to summation over $n=k+k' \geq 0$, we get:
\beq
\label{Q5}
\begin{array}{ll}
&F^{\zeta \xi}_{\zeta ' \xi '}(\lambda , \lambda ')=
\, \displaystyle{\lim _{\varepsilon \to 0}}\Bigl [ \varepsilon
{\cal X}^{\varepsilon}_{L}( \zeta +\lambda_{+}; \zeta , \xi ,
\lambda )\Bigr ] \, \times
\\&\\
\times &
\displaystyle{
\sum_{n\geq 0} \delta
(\zeta \!-\!\xi ' \!+\!\lambda \!+\!\lambda' \!-\!\eta \!-\!2n\eta )
\sum _{k=0}^{n} x_{k}(\zeta , \xi , \lambda )
{\cal X}_{R}( \zeta \!+\!\lambda _{+} \!-\!2k\eta ;
\zeta ' , \xi ' , \lambda ')} \,.
\end{array}
\eeq

The next step is to identify the finite sum over $k$
with the terminating modular hypergeometric function
${}_{10}W_{9}$. For the definition see Appendix A.
Plugging here the explicit formula for ${\cal X}_R$
(\ref{R2}) and using (\ref{gam6}), one can represent
(\ref{Q5}) in the form
\beq
\label{Q6}
\begin{array}{lll}
F^{\zeta \, \xi}_{\zeta ' \xi '} (\lambda , \lambda ')&=&
F^{\zeta \xi , \zeta ' \xi '}_{\mbox{{\small sym}}}
(\lambda , \lambda ')\,
\omega (\lambda ' , \zeta ' \!-\!\xi ')
\,F_{1}^{\zeta \xi}(\lambda)
F_{2}^{\zeta \zeta ' \xi '}(\lambda , \lambda ')
\,\displaystyle{
\frac{g_{\xi \xi '}(\lambda -\lambda ')}
{g_{\zeta \zeta '}(\lambda -\lambda ')} }\, \times
\\ && \\
&\times & \displaystyle{\sum_{n\geq 0} }
{}_{10}W_{9} (\alpha_1 ; \alpha_4, \ldots , \alpha_9, -n)
\delta
(\zeta \!-\!\xi ' \!+\!\lambda \!+\!\lambda' \!-\!\eta \!-\!2n\eta )\,,
\end{array}
\eeq
where
$$
F^{\zeta \xi , \zeta ' \xi '}_{\mbox{{\small sym}}}=
\frac{e^{2\pi i \zeta }\Gamma (2\zeta )}
{\Gamma (\zeta \!-\!\xi \!+\!2(\ell \!+\!1)\eta )
\Gamma (\zeta \!+\!\xi \!+\!2(\ell \!+\!1)\eta )}
\prod_{\epsilon =\pm }
\frac{\Gamma (\zeta +\epsilon \zeta '
+\lambda_{+}+\lambda_{+}' +\eta )}
{\Gamma (\zeta +\epsilon \xi '
+\lambda +\lambda ' +\eta )}
$$
is a symmetric function of $\lambda , \lambda '$ for all
$\zeta , \xi , \zeta ' , \xi '$,
$\omega (\lambda ' , \zeta ' \!-\!\xi ')$ is defined
in (\ref{omega}), and
\beq
\label{ggg}
g_{\zeta \zeta '}(\lambda )=\Gamma (\zeta +\zeta ' +\lambda +\eta )
\Gamma (\zeta -\zeta ' +\lambda +\eta )\,,
\eeq
\beq
\label{Q7}
F_{1}^{\zeta \xi}(\lambda )= e^{-4\pi i \ell \lambda }
\,\frac{\Gamma (\zeta -\xi +2\lambda ) \Gamma (\zeta +\xi +2\lambda )}
{\Gamma (2\lambda_{-})\Gamma (2\zeta +2\lambda _{+})}\,,
\eeq
\beq
\label{Q8}
F_{2}^{\zeta \zeta ' \xi '}(\lambda , \lambda ')=
\frac{\Gamma (\zeta +\xi ' +\lambda_{+}-\lambda_{-}' +\eta )
\Gamma (\zeta -\xi ' +\lambda_{+}-\lambda_{-}' +\eta )}
{\Gamma (\xi +\xi ' +\lambda -\lambda ' +\eta )
\Gamma (\xi -\xi ' +\lambda -\lambda ' +\eta )}\,.
\eeq
The function ${}_{10}W_{9}$ is defined in (\ref{mh3}).
It is an elliptic analog of the very-well poised hypergeometric
series. The values of parameters are as follows:
\beq
\label{Q9}
\begin{array}{ll}
-2\eta \alpha_1 =2\zeta \!+\!2\lambda \!+\!2\ell \eta&
-2\eta \alpha_7 =\zeta \!+\!\xi \!+\!2\ell \eta
\\
-2\eta \alpha_4 =2\lambda \!+\!2\ell \eta &
-2\eta \alpha_8 =\zeta \!-\!\xi \!+\!2\ell \eta
\\
-2\eta \alpha_5 =\zeta \!+\!\zeta ' \!+\!
\lambda \!-\!\lambda ' \!-\!\eta &
-2\eta \alpha_9 =\zeta \!+\!\xi ' \!+\!
\lambda \!+\!\lambda ' \!-\!\eta
\\
-2\eta \alpha_6 =\zeta \!-\!\zeta ' \!+\!
\lambda \!-\!\lambda ' \!-\!\eta &
-2\eta \alpha_{10} =\zeta \!-\!\xi ' \!+\!
\lambda \!+\!\lambda ' \!-\!\eta = 2n\eta
\end{array}
\eeq
Note that the parameters in the right column are symmetric
under the permutation $\lambda \leftrightarrow \lambda '$.
The series is {\it balanced} that means that the balancing
condition
\beq
\label{Q10}
2+3\alpha_1 =\sum_{m=1}^{7}\alpha_{m+3}
\eeq
is satisfied.
The series is terminating since $\alpha_{10}=-n$, where $n$
is a non-negative integer number.

Let us move the factor
$F_2 (\lambda , \lambda ')$ inside the sum $\sum_{n\geq 0}$
and set $\xi ' =\zeta +\lambda +\lambda ' -\eta -2n\eta$
in the $n$-th term (this is what
the $\delta$-function suggests to do). After using (\ref{gam6}),
(\ref{gam7}) and extracting the common factor
\beq
\label{Q11}
F_{3}^{\zeta \xi}(\lambda , \lambda ')=
\frac{\Gamma (2\zeta +2\lambda_{+})\Gamma (-2\lambda_{+}'+2\eta )}
{\Gamma (\zeta +\xi +2\lambda )
\Gamma (-\zeta +\xi -2\lambda '+2\eta )}\,,
\eeq
the sum over $n$ becomes equal to
\beq
\label{Q12}
F_{3}^{\zeta \xi}(\lambda , \lambda ')
\sum_{n\geq 0}
\frac{[\alpha_1 \!-\!\alpha_7\!+\!1]_n[\beta_1 \!-\!\alpha_7
\!-\!\alpha_8\!+\!1]_n}
{[\alpha_1 \!+\!1]_n[\beta_1 \!-\!\alpha_8\!+\!1]_n}
{}_{10}W_{9} (\alpha_1 ; \alpha_4, \ldots , -n)
\delta
(\zeta \!-\!\xi ' \!+\!\lambda \!+\!\lambda' \!-\!\eta \!-\!2n\eta )\,,
\eeq
where we use the elliptic Pochhammer symbol $[x]_n=[x][x+1]\ldots
[x+n-1]$, and $\beta_1=2\alpha_1\!+\!1\!-\!\alpha_4
\!-\!\alpha_5 \!-\!\alpha_6$.
Now, we recall the elliptic analog (\ref{Bailey})
of the Bailey transformation formula for the terminating
balanced ${}_{10}W_{9}$.
It is easily seen from (\ref{Q9}) that
the Bailey transformation of the
parameters $\alpha_i \to \beta_i$ (\ref{beta})
is equivalent to the permutation $\lambda \leftrightarrow \lambda '$
while all the other parameters remain unchanged.
Therefore, as it follows from (\ref{Bailey}),
each term inside the sum over $n$ in (\ref{Q12})
{\it is symmetric under the permutation
$\lambda \leftrightarrow \lambda '$}. The common factor
$F_3$ is not symmetric. Combining it with the remaining
factors and taking into account (\ref{gamma5}),
we obtain:
\beq
\label{Q13}
\begin{array}{ll}
&F^{\zeta \, \xi}_{\zeta ' \xi '} (\lambda , \lambda ')=
\displaystyle{
\, \frac{g_{\xi \xi '}(\lambda -\lambda ')}
{g_{\zeta \zeta '}(\lambda -\lambda ')} \, \times }\\&\\
\times & \displaystyle{\sum_{n\geq 0}}
{\cal F}^{\zeta \xi , \zeta ' \xi '}_{\mbox{{\small sym}}}
(\lambda , \lambda ';n)
\,\omega (\lambda , \zeta \!-\!\xi )\,
\omega (\lambda ' , \zeta '\!-\!\xi ')
\delta (\zeta \!-\!\xi ' \!+\!
\lambda \!+\!\lambda' \!-\!\eta \!-\!2n\eta )\,,
\end{array}
\eeq
where
${\cal F}^{\zeta \xi , \zeta ' \xi '}_{\mbox{{\small sym}}}
(\lambda , \lambda ' ;n)$
is symmetric under the permutation
of $\lambda$ and $\lambda '$ for all
$n, \zeta , \xi , \zeta ' , \xi '$, and the functions
$g$, $\omega$ are given in (\ref{ggg}), (\ref{omega})
respectively.

To obtain the kernel (\ref{Q2}), we should
set $\zeta =\zeta_i$, $\xi =\zeta_{i+1}$ and take
the product of the local factors (\ref{Q13})
over $i$ from $1$ to $N$ provided $\zeta_{N+1}=\zeta_1$.
Therefore, the ratios of the $g$-functions cancel.
The only factor which is not manifestly symmetric
under the permutation $\lambda \leftrightarrow \lambda '$
comes from the product of the $\omega$-functions.
This factor is equal to
\beq
\label{prom}
\exp \left ( -\, \frac{\pi i}{2\eta }
\sum_{j=1}^{N}(\zeta_{j}' -\zeta_{j+1}')^2 \right )
\prod_{i=1}^{N}
\left ( \theta_1(\zeta_i -\zeta_{i+1}-2\lambda |2\eta )\,
\theta_1(\zeta_i ' -\zeta_{i+1}'-2\lambda '|2\eta )\right )^{-1}\,.
\eeq
However, we should take into account the constraints
imposed by the $\delta$-functions:
$\zeta_i ' = \zeta_{i+1}-\lambda -\lambda ' +\eta +2n_{i-1}\eta$,
where $n_i \in \ZZ$. In other words, (\ref{prom}) is zero unless
$\zeta_i ' -\zeta_{i+1}'
= \zeta_{i-1}-\zeta_i +2(n_{i-1}-n_i )\eta$.
From the monodromy properties of the $\theta$-function
it then follows that this function is symmetric, too.
Therefore, the commutation relation (\ref{Q1})
for ${\bf Q}_{L,R}^{+}$ is proved. To prove it for
${\bf Q}_{L,R}^{-}$, there is no need to repeat the
calculations since $X_{L,R}^{-}$ is related to
$X_{L,R}^{+}$ via (\ref{E2}).

The proof of (\ref{Q1a}) requires a similar but different
calculation. Its scheme is outlined below.
The kernel of the operator in the left hand side is
\beq
\label{Q2a}
\begin{array}{ll}
&\displaystyle{\prod_{i=1}^{N}
\int dz_i X_L^+ (\lambda , \zeta_i , \zeta_{i+1};z_i )
X_R^- (z_i ; \zeta'_i , \zeta'_{i+1}, \lambda ')}\,=
\\&\\
=&\displaystyle{\prod_{i=1}^{N}
\int dz_i X_L^+ (\lambda , \zeta_i , \zeta_{i+1};z_i )
X_R^+ (z_i ; \zeta'_{i+1}, \zeta'_{i}, -\lambda ')}\,,
\end{array}
\eeq
where (\ref{E2}) is used.
Again, set $\zeta_i =\zeta$, $\zeta_{i+1}=\xi$,
$\zeta_i' =\zeta '$, $\zeta_{i+1}'=\xi '$,
and examine the symmetry properties of the function
$$
\tilde F^{\zeta \xi}_{\zeta ' \xi '}(\lambda , \lambda ')=
\int dz X_L^+ (\lambda , \zeta , \xi ;z )
X_R^+ (z ; \xi , \zeta , -\lambda ')
$$
under the simultaneous change $\zeta \leftrightarrow \xi$,
$\zeta ' \leftrightarrow \xi '$, $\lambda \to -\lambda '$,
$\lambda ' \to -\lambda$. It is enough to show that
the product of such functions over the chain is invariant
under this substitution in each factor.
Proceeding as above, we obtain the representation
of this function similar to (\ref{Q6}):
\beq
\label{Q6a}
\begin{array}{lll}
\tilde F^{\zeta \, \xi}_{\zeta ' \xi '} (\lambda , \lambda ')&=&
F^{\zeta \xi , \zeta ' \xi '}_{\mbox{{\small sym}}}
(\lambda , \lambda ')\,
\omega (-\lambda ' , \xi ' \!-\!\zeta ')
\,F_{1}^{\zeta \xi}(\lambda)
F_{2}^{\zeta \zeta ' \xi '}(\lambda , \lambda ')
\,\displaystyle{
\frac{g_{\xi \xi '}(\lambda -\lambda ')}
{g_{\zeta \zeta '}(\lambda -\lambda ')} }\, \times
\\ && \\
&\times & \displaystyle{\sum_{n\geq 0} }
{}_{10}W_{9} (\tilde \alpha_1 ; \tilde \alpha_4, \ldots ,
\tilde \alpha_{10})
\delta
(\zeta \!-\!\zeta ' \!+\!\lambda \!-\!\lambda' \!-\!\eta \!-\!2n\eta )
\end{array}
\eeq
with the same functions $F_{\mbox{{\small sym}}}$,
$\omega$, $F_1$, $F_2$, and $g$ as before.
The parameters of the series ${}_{10}W_{9}$
are also the same (see (\ref{Q9}));
however, for the sake of convenience they are permuted
in the following way: $\tilde \alpha_1 =\alpha_1$,
$\tilde \alpha_4 =\alpha_8$,
$\tilde \alpha_5 =\alpha_9$,
$\tilde \alpha_6 =\alpha_{10}$,
$\tilde \alpha_7 =\alpha_4$,
$\tilde \alpha_8 =\alpha_5$,
$\tilde \alpha_9 =\alpha_7$,
$\tilde \alpha_{10} =\alpha_6$. (Clearly, the series is
symmetric under the permutation of all parameters except
$\alpha_1$.)
What differs is the
argument of the $\delta$-function.
Note that the series is terminating since
the $\delta$-function is zero unless $\tilde \alpha_{10} =-n$.
Next, we make the Bailey transformation (\ref{Bailey})
passing from $\tilde \alpha_i$ to $\tilde \beta_i$
according to (\ref{beta}). The result can be represented
as a product of three factors:
$$
\tilde F^{\zeta \, \xi}_{\zeta ' \xi '} (\lambda , \lambda ')
= \,A^{\zeta \, \xi}_{\zeta ' \xi '} (\lambda , \lambda ')\,\times \,
\left [ \begin{array}{l}\mbox{invariant under}
\,\,\zeta \!\leftrightarrow \!\xi, \\
\zeta ' \! \leftrightarrow \! \xi ',\,
\lambda \!\leftrightarrow \!- \lambda '
\end{array} \right ] \times
\left [ \begin{array}{l}
\mbox{factors that cancel in the}\\
\mbox{product over the chain}
\end{array} \right ],
$$
where
$$
A^{\zeta \, \xi}_{\zeta ' \xi '} (\lambda , \lambda ')=
\frac{\exp ( -\pi i (\xi ' -\zeta ' )^2 /(2\eta ))}
{\theta_1 (\xi ' \!-\! \zeta ' \!-\! 2\lambda ' |2\eta )
\theta_1 (\zeta \!-\! \xi +\! 2\lambda |2\eta )}\,.
$$
The product $\prod_{i=1}^{N}
A^{\zeta_i \, \zeta_{i+1}}_{\zeta_i' \zeta_{i+1}'}
(\lambda , \lambda ')$ is the only thing here which is not
manifestly invariant
under the simultaneous change
$\zeta_{i}\leftrightarrow \zeta_{i+1}$,
$\zeta_{i}'\leftrightarrow \zeta_{i+1}'$,
$\lambda \leftrightarrow -\lambda '$ in each factor.
However, the same argument as the one
used for (\ref{prom}) establishes the invariance.

\section{From pre-$Q$-operators to $Q$-operator}

Suppose the operators ${\bf Q}^{+}_L$
and ${\bf Q}^{+}_R$
are invertible at some
particular $\lambda =\lambda_{*}$. Following \cite{Baxter},
we define
\beq
\label{N1}
{\bf Q}^{\pm}(\lambda )=({\bf Q}^+_{L}(\lambda_{*}))^{-1}
{\bf Q}^{\pm}_{L}(\lambda ) =
{\bf Q}^{\pm}_{R}(\lambda )({\bf Q}^+_{R}(\lambda_{*}))^{-1}\,.
\eeq
The second equality is due to the basic commutation
relations (\ref{Q1}), (\ref{Q1a}).
It then follows, on applying (\ref{Q1}), (\ref{Q1a})
once again, that $[{\bf Q}^{+}(\lambda ),\,{\bf Q}^{+}(\mu )]
=[{\bf Q}^{+}(\lambda ),\,{\bf Q}^{-}(\mu )]
=[{\bf Q}^{-}(\lambda ),\,{\bf Q}^{-}(\mu )]=0$
for all $\lambda$, $\mu$:
\beq
\label{N2}
\begin{array}{ll}
&{\bf Q}^{+}(\lambda ){\bf Q}^{+}(\mu )\, =
\\&\\
=& ({\bf Q}^{+}_{L}(\lambda_{*}))^{-1}
{\bf Q}^{+}_{L}(\lambda ) {\bf Q}^{+}_{R}(\mu )
({\bf Q}^{+}_{R}(\lambda_{*}))^{-1}\, =
\\&\\
=& ({\bf Q}^{+}_{L}(\lambda_{*}))^{-1}
{\bf Q}^{+}_{L}(\mu ) {\bf Q}^{+}_{R}(\lambda )
({\bf Q}^{+}_{R}(\lambda_{*}))^{-1}\, =
\\&\\
=&{\bf Q}^{+}(\mu ){\bf Q}^{+}(\lambda )\,,
\end{array}
\eeq
and similarly for ${\bf Q}^{-}$.
Both equations (\ref{TQL}), (\ref{TQR}) lead to the
Baxter $T$-$Q$ relation
\beq
\label{TQ}
{\bf T}(\lambda ){\bf Q}(\lambda )=
\theta_{1}^{N}(2\lambda -2\ell \eta )
{\bf Q}(\lambda +\eta )
+\theta_{1}^{N}(2\lambda +2\ell \eta )
{\bf Q}(\lambda -\eta )\,,
\eeq
where ${\bf Q}$ is any one of ${\bf Q}^{\pm}$.
Besides, it follows that ${\bf Q}^{\pm}(\lambda )$ commute
with ${\bf T}(\mu )$ for all $\lambda$, $\mu$.

Suppose $\ell \notin \ZZ_{+}$. Then the convenient
choice of the normalization point $\lambda_{*}$
is $\lambda_{*}=-\ell \eta$. Indeed,
at $\lambda =-\ell \eta$ the series (\ref{Q4})
terminates at the $k=0$ term; therefore,
${\bf Q}^{+}(-\ell \eta )$ is a diagonal operator.
Since we did not take care of $\ell$-dependent constant
factors, we continue to do so and redefine the $Q$-operator
dividing it by the constant $\Gamma^{N}(-4\ell \eta )$.
Note that this normalization allows us to define the
$Q$-operator for $\ell \in \ZZ_{+}$ as well.
We express the result through the function $\Phi_{\ell}$
(\ref{Phig}):
\beq
{\bf Q}_{L}^{+}(-\ell \eta ) =
%\Gamma^{-N}(-4\ell \eta )
\prod_{j=1}^{N}
\Phi_{\ell}^{-1}(z_j -z_{j+1})
\Phi_{\ell}^{-1}(z_j +z_{j+1})\,.
%\frac{e^{2\pi i (2\ell +1)z_j}
%\Gamma (z_j \!-\!z_{j+1}\!-\!2\ell \eta )
%\Gamma (z_j \!+\!z_{j+1}\!-\!2\ell \eta )}
%{\Gamma (z_j -z_{j+1}+2(\ell +1)\eta )
%\Gamma (z_j +z_{j+1}+2(\ell +1)\eta )}\,.
\eeq

Thus, by virtue of (\ref{kQL}), the kernels
of the $Q$-operators ${\bf Q}^{\pm}$ are
\beq
\label{N3}
Q^{\pm}(\lambda ; \vec \zeta , \vec z)=
\prod_{i=1}^{N}q^{\pm}(\lambda ; \zeta_i , \zeta_{i+1},z_i )\,,
\eeq
where
\beq
\label{N4}
\begin{array}{lll}
q^{+}(\lambda ; \zeta_i , \zeta_{i+1},z_i )&=&
\displaystyle{
\frac{\Gamma (2\zeta_i )
\Gamma (\zeta_i -\zeta_{i+1}+2\lambda )
\Gamma (\zeta_i +\zeta_{i+1}+2\lambda )}
{\Gamma (2\lambda_{-}) \Gamma (2\zeta_i \!+\!2\lambda_{+})
\Gamma (\zeta_i \!-\!\zeta_{i+1}\!-\!2\ell \eta )
\Gamma (\zeta_i \!+\!\zeta_{i+1}\!-\!2\ell \eta )} }\,\,\times
\\&& \\
&\times & \displaystyle{
\sum_{k\geq 0}x_k (\zeta_i , \zeta_{i+1}, \lambda )}
\delta (z_i -\zeta_i -\lambda_{+} +2k\eta )
\end{array}
\eeq
($x_k$ is defined in (\ref{xk})), and
$q^{-}(\lambda , \zeta_i , \zeta_{i+1}, z_i ) =
q^{+}(-\lambda , \zeta_{i+1}, \zeta_{i}, z_i)$.
Representing this in the operator form, according to
the rule described in Sec.\,3, we obtain:
\beq
\label{N6}
\begin{array}{lll}
{\bf Q}^{+}(\lambda )&=&\!\!
\displaystyle{
\frac{1}{\Gamma^N(2\lambda \!-\!2\ell \eta )}
\prod_{i=1}^{N}
\frac{\Gamma (2z_i )\,\Gamma (z_i \!-\!z_{i+1}\!+\!2\lambda )
\Gamma (z_i \!+\!z_{i+1}\!+\!2\lambda )}
{\Gamma (2z_i \!+\!2\lambda \!+\!2\ell \eta )
\Gamma (z_i \!-\!z_{i+1}\!-\!2\ell \eta )
\Gamma (z_i \!+\!z_{i+1}\!-\!2\ell \eta )}}
\,\, \times \\ && \\
&\times & \!\! \!\!\normord  \,
\!\displaystyle{\prod_{j=1}^{N}}
{}_{6}W_{5} \Bigl (   \displaystyle{
-\frac{z_j\!+\!\lambda }{\eta}\!-\!\ell ;\,
-\frac{\lambda }{\eta}\!-\!\ell ,\,
-\frac{z_{j+1}\!+\!z_j}{2\eta}\!-\!\ell ,\,
\frac{z_{j+1}\!-\!z_j}{2\eta}\!-\!\ell ; \,
e^{-2\eta \p _{z_j}} }\Bigl | \eta , \tau \Bigr ) \normord  \,\times
\\&&\\
&\times & \exp \Bigl ((\lambda +\ell \eta )
\displaystyle{\sum_{m=1}^{N}\p_{z_m} }\Bigr )\,,
\end{array}
\eeq
\beq
\begin{array}{lll}
\label{N7}
{\bf Q}^{-}(\lambda )&=&\!\!
\displaystyle{
\frac{1}{\Gamma^N(-2\lambda \!-\!2\ell \eta )}
\prod_{i=1}^{N}
\frac{\Gamma (2z_i )\,\Gamma (z_{i} \!-\!z_{i-1}\!-\!2\lambda )
\Gamma (z_{i} \!+\!z_{i-1}\!-\!2\lambda )}
{\Gamma (2z_i \!-\!2\lambda \!+\!2\ell \eta )
\Gamma (z_{i}\!-\!z_{i-1}\!-\!2\ell \eta )
\Gamma (z_{i} \!+\!z_{i-1}\!-\!2\ell \eta )}}
\,\, \times \\ && \\
&\times & \!\! \!\!\normord  \,
\!\displaystyle{\prod_{j=1}^{N}}
{}_{6}W_{5} \Bigl (   \displaystyle{
-\frac{z_{j}\!-\!\lambda }{\eta}\!-\!\ell ;\,
\frac{\lambda }{\eta}\!-\!\ell ,\,
-\frac{z_{j}\!+\!z_{j-1}}{2\eta}\!-\!\ell ,\,
\frac{z_{j}\!-\!z_{j-1}}{2\eta}\!-\!\ell ; \,
e^{-2\eta \p _{z_{j}}} }\Bigl | \eta , \tau \Bigr ) \normord  \,\times
\\&&\\
&\times & \exp \Bigl ((-\lambda +\ell \eta )
\displaystyle{\sum_{m=1}^{N}\p_{z_m} }\Bigr )\,{\bf C}\,,
\end{array}
\eeq
where ${\bf C}$ is the cyclic shift operator,
\beq
\label{C}
({\bf C}f)(z_1 , z_2 , \ldots , z_{N-1}, z_{N})=
f(z_2 , z_3 , \ldots , z_{N}, z_{1})\,,
\eeq
with the kernel
$C(\vec \zeta , \vec z)=\prod_{i=1}^{N} \delta (z_i -\zeta_{i+1})$.
Clearly, ${\bf C}$ commutes with the ${\bf T}(\lambda )$ and
the ${\bf Q}^{\pm}(\lambda )$.

In general, the series
${}_{6}W_{5}$ in
(\ref{N6}), (\ref{N7}) is not terminating. However,
at $\lambda =\mp \ell \eta \pm m\eta$, $m\in \ZZ_{+}$, the
series ${}_{6}W_{5}$ in (\ref{N6})
(respectively, in (\ref{N7})) terminates at the $m$-th
term. Note also that
\beq
{\bf Q}^{\pm}(\pm \ell \eta \mp m\eta )=0\,,
\;\;\;\;\; m=0, 1, 2 ,\ldots \,,
\label{zeros}
\eeq
due to zeros of the scalar common factor. At the points
$\pm \ell \eta$ we have:
\beq
\label{0011}
{\bf Q}^{+}(-\ell \eta )=
\frac{1}{\Gamma ^N(-4\ell \eta )}\,{\bf I}\,,
\;\;\;\;\;
{\bf Q}^{-}(\ell \eta )=
\frac{1}{\Gamma ^N(-4\ell \eta )}\,{\bf C}\,,
\eeq
where ${\bf I}$ is the identity operator. So, at
$\ell \in \frac{1}{2}\ZZ_{+}$ the $Q$-operators have poles
at these points.

Let us evaluate the $Q$-operator at another
distinguished point.
Note that if $\lambda =(\ell \!+\!1)\eta $, then
$\zeta_i$ in the kernel of ${\bf Q}^{+}$
decouples from $\zeta_{i+1}$, so
the operator simplifies:
\beq
{\bf Q}^{+}((\ell +1)\eta )=
\Bigl (\prod_{i=1}^{N}
\Phi_{\ell}(z_i -z_{i+1})\Phi_{\ell}(z_i +z_{i+1})
\Bigr )
\prod_{j=1}^{N}\stackrel{j}{ {\bf w}}_{\ell}\,.
\eeq
The operators ${\bf w}_{\ell}$
pairwise commute since they act on variables attached
to different sites of the chain.
From (\ref{N6}) one obtains
the representation of the operator ${\bf w}_{\ell}$
through the elliptic hypergeometric series
${}_{4}W_{3}$ with an operator argument:
\beq
\label{winf}
{\bf w}_{\ell}=
\frac{e^{2\pi i (2\ell +1)z}\Gamma (2z)}{\Gamma (2z
+2(2\ell +1)\eta )}
\normord
{}_{4}W_{3}\Bigl (-\frac{z}{\eta}\!-\!2\ell \!-\!1 ;\,
-\!2\ell \!-\!1 ; \,e^{-2\eta \p_{z}}\Bigl | \eta , \tau \Bigr )
\normord
e^{(2\ell +1)\eta \p_{z}} \,.
\eeq
At $\ell \in \frac{1}{2}\ZZ _{+}$ the series is terminating
and gives the operator ${\bf w}_{\ell}$ introduced
in (\ref{W}).
So, at $\lambda =(\ell +1)\eta$ we obtain a distinguished
factorized integral of motion
\beq
\label{A}
{\bf A}=\Bigl ( \prod_{i=1}^{N}\Phi_{\ell}(z_i +z_{i+1})
\Phi_{\ell}(z_i -z_{i+1})\Bigr )
\prod_{j=1}^{N}\stackrel{j}{{\bf w}}_{\ell}
\eeq
(no normal ordering is required
since the $\stackrel{\!j}{{\bf w}_{\ell}}$ commute).

Let us give some remarks on the operator ${\bf A}$.
It is clear from Sec.\,3 that for $\ell \in \frac{1}{2}\ZZ_{+}$
the operator ${\bf A}$ annihilates the space
$(\Theta_{4\ell}^{+})^{\otimes N}$.
At $\ell =0$, this integral of motion
is proportional to ${\bf T}(0)$ and has
especially simple form:
$2^N\Bigl ( \prod_{i=1}^{N}\theta_1(z_i \!+\!z_{i+1})
\theta_1(z_i \!-\!z_{i+1})/\theta_1(2z_i)\Bigr )
\prod_{j=1}^{N}\,\mbox{sinh}\,(\eta \p_z )$.
(In the limit to the $XXX$ spin chain it coincides with
the integral of motion found by Li\-pa\-tov \cite{Lipatov}
in the context of high energy QCD, see also \cite{FK}.) Presumably,
for $\ell \in \frac{1}{2}\ZZ_{+}$ ${\bf A}$
can be represented as the trace
of the product of more general ${\sf L}$-operators
obtained from (\ref{L}) by the fusion procedure, taken
at a particular value of $\lambda$.
For brevity, we denote
$\varphi_{\ell}(z_i , z_{i+1})=
\Phi_{\ell}(z_i \!-\!z_{i+1})\Phi_{\ell}(z_i \!+\!z_{i+1})$.
Together with (\ref{S3}),
the commutativity of
${\bf A}$ and ${\bf T}(\lambda )$ implies the intertwining
relation
\beq
\label{T5}
{\bf T}^{(\ell )}(\lambda)\,\prod_{i=1}^{N}\varphi_{\ell}(z_i , z_{i+1})
=\Bigl (\prod_{i=1}^{N}\varphi_{\ell}(z_i , z_{i+1})\Bigr )
{\bf T}^{(-\ell -1)}(\lambda)\,,
\eeq
where the superscript of ${\bf T}$ indicates the value of spin.
A direct proof of this relation
is given in Appendix B. We stress that
the intertwining relation (\ref{S3}) holds for all the spin
operators separately while (\ref{T5}) is valid for the
transfer matrix only.

The operators ${\bf Q}^{\pm}(\lambda )$ are two
operator solutions of the Baxter $T$-$Q$ relation (\ref{TQ})
treated as a difference equation in $\lambda$. Therefore,
their Wronskian,
$$
{\bf W}(\lambda ) ={\bf Q}^{+}(\lambda +\eta )
{\bf Q}^{-}(\lambda )
-{\bf Q}^{+}(\lambda )
{\bf Q}^{-}(\lambda +\eta )\,,
$$
obeys the functional equation
$\theta_{1}^{N}(2\lambda -2\ell \eta ){\bf W}(\lambda )=
\theta_{1}^{N}(2\lambda +2\ell \eta ){\bf W}(\lambda -\eta )$.
Taking into account (\ref{zeros}),
one can write
\beq
\label{wr}
{\bf W}(\lambda )=
\theta_1^N(2\lambda \!+\! 2\ell \eta |2\eta )\left (
\frac{\Gamma (2\lambda +2\ell \eta +2\eta )}{\Gamma (2\lambda
-2\ell \eta +2\eta )}
\right )^N {\bf p}(\lambda ){\bf A C}\,,
\eeq
where ${\bf p}(\lambda )$ is a diagonal operator
such that ${\bf p}(\lambda +1)={\bf p}(\lambda )$,
${\bf p}(\lambda +\eta )=e^{2\pi i N \eta}{\bf p}(\lambda )$, and
${\bf p}(\ell \eta )={\bf I}$.
The operator ${\bf A}$ is given in (\ref{A}).
(The simplest way to find the operator factors is to calculate
the Wronskian at $\lambda =\ell \eta$.)
To find ${\bf p}(\lambda )$ explicitly, a more careful analysis
is required.

\section{Miscellanea}

\subsubsection*{Twisted boundary conditions}

The twisted boundary condition means inserting
a $c$-number matrix inside the trace of the chain of
${\sf L}$-operators (\ref{T3}).
To preserve integrability, this matrix must be
a constant solution to the
${\sf R}{\sf L}{\sf L}={\sf L}{\sf L}{\sf R}$
relation. In the elliptic case,
there are only three such solutions: the Pauli
matrices $\sigma_{\alpha}$, $\alpha =1,2,3$.
Therefore, the definition of the quantum transfer
matrix with twisted boundary conditions is as follows:
\beq
\label{tw1}
{\bf T}(\lambda )=\mbox{tr}\, \bigl [
\stackrel{1}{{\sf L}}(\lambda)
\stackrel{2}{{\sf L}}(\lambda) \, \ldots \,
\stackrel{N}{{\sf L}}(\lambda) \sigma_{\alpha}\bigr ]\,.
\eeq
For definitness, we consider the case of
$\sigma_1= \left ( \begin{array}{cc}0&1 \\ 1&0 \end{array}
\right )$.

Again, we employ the gauge transformation (\ref{G1})
with the same matrix (\ref{G2}). The whole argument remains
the same with the only modification that $\zeta_{N+1}\neq \zeta_1$.
The $\zeta_{N+1}$
should be adjusted in such a way that the
matrix $\sigma_{1}'= G(\zeta_{N+1})\sigma_1
(G(\zeta_1 ))^{-1}$ be lower-triangular.
Calculating it explicitly, we obtain
\beq
\sigma_{1}'= \frac{1}{\bar \theta_4 (\zeta_1 )}
\left ( \begin{array}{cc}
\bar \theta_3 (\zeta_{N+1}) &
-2\theta_2 (\zeta_1 \!+\! \zeta_{N\!+\!1})
\theta_2 (\zeta_1 \!-\! \zeta_{N\!+\!1}) \\ &\\
1& -\bar \theta_3 (\zeta_{1}) \end{array} \right ),
\eeq
where a simple identity for $\theta$-functions has been used.
Therefore, if $\zeta_{N+1}=\zeta_1 +\frac{1}{2}$, then
the right upper element is zero.

With this modification, the definitions
of the vacuum vectors (\ref{G3R}), (\ref{G3L}) and
of the pre-$Q$-operators
(\ref{kQL}), (\ref{kQR}) remain the same.
However, equations (\ref{TQL}), (\ref{TQR}) get modified:
\beq
\label{TQLa}
{\bf T}(\lambda) {\bf Q}^{\pm}_{R}(\lambda)=
\theta_{1}^{N}(2\lambda_{-} )
{\bf Q}^{\pm}_{R}(\lambda +\eta)-
\theta_{1}^{N}(2\lambda_{+})
{\bf Q}^{\pm}_{R}(\lambda -\eta)\,,
\eeq
\beq
\label{TQRa}
{\bf Q}^{\pm}_{L}(\lambda) {\bf T}(\lambda)=
\theta_{1}^{N}(2\lambda_{-} )
{\bf Q}^{\pm}_{L}(\lambda +\eta)-
\theta_{1}^{N}(2\lambda_{+})
{\bf Q}^{\pm}_{L}(\lambda -\eta)\,.
\eeq

The new constraint on $\zeta_{N+1}$ and $\zeta_1$
does not destroy the basic commutation
relation (\ref{Q1}). Indeed, the $g$-factors in (\ref{Q13})
still cancel in the product
over the chain since $g_{\zeta_1 \zeta_{1}'}=
g_{\zeta_{N+1} \zeta_{N+1}'}$ due to (\ref{ggg}) and
(\ref{gamma1}).

The normalization of the $Q$-operators goes through as before.
The $T$-$Q$ relation is
\beq
\label{TQa}
{\bf T}(\lambda ){\bf Q}(\lambda )=
\theta_{1}^{N}(2\lambda -2\ell \eta )
{\bf Q}(\lambda +\eta )
-\theta_{1}^{N}(2\lambda +2\ell \eta )
{\bf Q}(\lambda -\eta )\,.
\eeq
The $Q$-operators are given by the same formulas
(\ref{N6}) (\ref{N7}),
where it is implied that $\zeta_{N+1}=\zeta_1 +\frac{1}{2}$,
and the cyclic shift operator ${\bf C}$ should be replaced
by the operator ${\bf C}'$ such that
$({\bf C}'f)(z_1 , z_2 , \ldots , z_{N})=
f(z_2 , z_3 , \ldots , z_{N}, z_{1}+\frac{1}{2})$.

\subsubsection*{The case $N=1$}

For $N=1$, there is only one ${\sf L}$-operator
in the chain (\ref{T3}).
Trace of the $\mbox{{\sf L}}(\lambda)$
is the simplest
quantum transfer matrix ${\bf T}(\lambda)$. It is proportional
to ${\bf s}_0$:
\beq
\label{T1}
{\bf T}(\lambda)=\mbox{tr}\,\mbox{{\sf L}}(\lambda)
=\theta _{1}(2\lambda){\bf s}_0\,,
\eeq
\beq
{\bf s}_0 = \frac{\theta _{1}(2z -2\ell \eta)}
{\theta _{1}(2z)}\, e^{\eta \p _z}
+\frac{\theta _{1}(2z +2\ell \eta)}
{\theta _{1}(2z)}\, e^{-\eta \p _z}\,.
\label{Int1}
\eeq
This is the
difference analog of the Lam\'e operator
proposed in \cite{kz}.
The difference operator given by eq.\,(\ref{Int1})
for positive integer values
of $\ell$ and {\it arbitrary} generic $\eta$
has $2\ell +1$ stable bands (and $2\ell$ gaps) in the spectrum.

The construction of the previous sections gives
a remarkable family
of commuting difference operators in one variable.
For $\ell \in \ZZ_{+}$ this family of operators commuting with
(\ref{Int1}) coincides with the one found in \cite{FV2}
by means of a different argument.

The case $N=1$ is exceptional and must be considered
separately. Some of the
arguments for $N\geq 2$ are not directly applicable to
$N=1$ since the vacuum vectors (\ref{R3}), (\ref{L3})
are to be taken at $\zeta =\xi$ that for
$\ell \in \ZZ_{+}$ means the
{\it non-generic position} of the parameters.
Nevertheless, it is a matter of direct verification
to prove that
the final formulas (\ref{N6}), (\ref{N7}) remain true,
if just to change the normalization making them
non-singular for $\ell \in \ZZ_{+}$:
\beq
\begin{array}{lll}
{\bf Q}^{\pm}_{(N=1)}(\pm \lambda )&=&\!\!
\displaystyle{
\frac{\Gamma (2\lambda )\,\Gamma (2z) \,\Gamma (2z+2\lambda )}
{\Gamma (2\lambda \!-\! 2\ell \eta )
\Gamma (2z \!-\! 2\ell \eta )
\Gamma (2z \!+\!2\lambda \!+\! 2\ell \eta )} }
\,\, \times \\ && \\
&\times & \!\! \normord  \,
{}_{6}W_{5} \Bigl (   \displaystyle{
-\frac{z+\lambda }{\eta}-\ell ;\,
-\frac{\lambda }{\eta}-\ell ,\,
-\frac{z}{\eta}-\ell ,\, -\ell ; \,
e^{-2\eta \p _{z}} }\Bigl | \eta , \tau \Bigr ) \normord  \,
e^{(\lambda +\ell \eta )\p_z}\,.
\end{array}
\eeq
For $\ell \in \ZZ _{+}$ the series is terminating. In this case
the above formula can be rewritten in a more explicit form:
\beq
\label{op1}
{\bf Q}^{\pm}_{(N=1)}(\lambda )=
\sum _{k=0}^{\ell}
A_{k}(z,\lambda)
\,e^{(2k\eta -\ell \eta +\lambda)\p _z}\,,
\eeq
where
\begin{eqnarray}
\label{op2}
\!\!\!A_{k}(z,\lambda)\!\!
&\,=\,(-1)^k
\displaystyle{\frac{[\ell ]!}{[2\ell ]!}}
\left [\begin{array}{c}\ell \\ k\end{array} \right ]
&\,\prod _{j=0}^{\ell -k-1}
\frac{\theta _{1}\bigl (2z +2(\ell \!-\!j)\eta \bigr )
\theta _{1}\bigl (2\lambda +2(\ell \!-\!j)\eta \bigr )}
{\theta _{1}\bigl (2z +2\lambda +2(k \!-\!j)\eta \bigr )}
\nonumber\\
&&\!\!\!\times\,\,\prod _{j=0}^{k-1}
\frac{\theta _{1}\bigl (2z -2(\ell \!-\!j)\eta \bigr )
\theta _{1}\bigl (2\lambda -2(\ell \!-\!j)\eta \bigr )}
{\theta _{1}\bigl (2z +2\lambda +2(k \!+\!j\!-\!\ell )
\eta \bigr )}\,.
\end{eqnarray}
If $k=0$ or $k=\ell$, then the second (respectively, the first)
product is absent.

\subsubsection*{The $XXZ$-case}

The construction of Secs.\,5 -- 7 admits
several different trigonometric degenerations. Here we discuss
the one corresponding to the $XXZ$ spin chain with the standard
6-vertex $R$-matrix. The ${\sf L}$-operator is (cf.\,(\ref{L}))
\beq
\label{Z1}
\mbox{{\sf L}}^{(XXZ)}(u)=\frac{1}{2i}
\left ( \begin{array}{cc}
e^{2\pi i u+{\bf j}_0}
-e^{-2\pi i u-{\bf j}_0}
& {\bf j}_{-}
\\& \\
{\bf j}_{+}
&
e^{2\pi i u-{\bf j}_0}
-e^{-2\pi i u +{\bf j}_0}
\end{array} \right ),
\eeq
where $u$ is the spectral parameter and
${\bf j}_0$, ${\bf j}_{\pm}$ are
operators acting on functions of a variable $x$ as follows:
\beq
\label{Z2}
{\bf j_0}=\gamma \p_x\,, \;\;\;\; {\bf j}_{\pm} =e^{\pm 2\pi i x}
(q^{\ell}e^{\mp \gamma \p_x} -q^{-\ell}e^{\pm \gamma \p_x}) \,.
\eeq
Here $q=e^{2\pi i \gamma}$,
and these operators obey the commutation relations of
the quantum algebra $U_q(sl(2))$:
\beq \label{Z3} [{\bf j}_0 ,\, {\bf j}_{\pm}]=
\pm 2\pi i\gamma {\bf j}_{\pm}\,,
\;\;\;\;[{\bf j}_{+} ,\, {\bf j}_{-}]= 2(q-q^{-1})\,
\mbox{sinh}\, {\bf j}_0 \,.
\eeq
The quantum transfer matrix ${\bf T}(u)$ is defined by the same
formula (\ref{T3}).

The limiting transition from the elliptic formulas to the
trigonometric ones consists of two steps:
\begin{itemize}
\item[i)]
The modular transformation $\tau \to -1/\tau$ with the change
of variables $z_j=\tau x_j +\frac{1}{4}$, $\lambda =\tau u$,
$\eta =\tau \gamma$;
\item[ii)]
Tending $\tau$ to $0$.
\end{itemize}
Then the ${\sf L}$-operator defined as ${\sf L}^{(XXZ)}(u)=
\lim_{\tau \to i0}\sqrt{\frac{i\tau}{4}}\,e^{\frac{i\pi}{4\tau}}
{\sf L}(\lambda )$ with the understanding that
$\tau$ in the ${\sf L}(\lambda )$ (\ref{L}) is replaced by
$-1/\tau$, and
$u=\lambda /\tau$
is kept fixed coincides with (\ref{Z1}).
The same limiting transition in (\ref{N6}) yields the
trigonometric $Q$-operator.

Skipping the details (for the modular transformation
of the $\theta$-functions and the elliptic gamma-function
see eqs.\,(\ref{mod}), (\ref{modg}), respectively), we
present the final result. It is written in terms of the
following special functions. The basic hypergeometric
function ${}_{2}\phi _{1}$ \cite{GR} is defined as
\beq
\label{basic}
{}_{2}\phi _{1}\left [\begin{array}{c}a,b\\c \end{array};
q, x\right ]=
\sum_{k=0}^{\infty} \frac{(a;q)_k (b;q)_k}
{(q;q)_k (c;q)_k}\, x^k\,,
\eeq
where $(a;q)_k =(1-a)(1-aq)(1-aq^2)\ldots (1-aq^{k-1})$.
Up to simple redefinitions of parameters, this series
gives the little $q$-Jacobi function (when the series
is terminating, it is the little $q$-Jacobi polynomial).
The $q$-gamma-function is
\beq
\label{q-gam}
\Gamma_q(x)=\frac{(q;q)_{\infty}}{(q^x;q)_{\infty}}\,
(1-q)^{1-x}\,,
\eeq
where, for $|q|<1$,
$(a;q)_{\infty}=\lim_{k\to \infty}(a;q)_k$.
The $Q$-operator is given by
\beq
\label{Z4}
\begin{array}{lll}
{\bf Q}^{+}(u)&=&e^{-2\pi i N\ell u}\,
\displaystyle{\prod_{i=1}^{N}}
\displaystyle{
\frac{\Gamma_{q^2}
\left (
\frac{ x_i -x_{i+1}+2u}{2\gamma}\right )}
{\Gamma_{q^2} \left ( \frac{u}{\gamma} \!-\!\ell \right )
\Gamma_{q^2}
\left ( \frac{x_i - x_{i+1} }{2\gamma}\!-\!\ell \right )} }
\, \times \\&&\\
&\times & \normord
\displaystyle{\prod_{j=1}^{N} }
{}_{2}\phi _{1}\left [
\begin{array}{c}
q^{-2\ell}e^{-4\pi i u},\,
q^{-2\ell}e^{2\pi i (x_{j+1}\!-\!x_j )}\\
q^{2}e^{2\pi i (x_{j+1}\!-\!x_j \!-\!2u)}
\end{array}; q^2 ,\, q^{2(\ell +1)}e^{-2\gamma \p_{x_j}} \right ]
\normord \, \times
\\&&\\
&\times &\exp \Bigl ((u+\ell \gamma )\sum_{m=1}^{N}\p_{x_m}\Bigr )\,.
\end{array}
\eeq
As in (\ref{N7}), the formula for ${\bf Q}^{-}(u)$ differs
from (\ref{Z4}) by the changes $u\to -u$,
$x_i \leftrightarrow x_{i+1}$, and the additional cyclic shift
${\bf C}$. The both operators ${\bf Q}^{\pm}(u)$
obey the $T$-$Q$ relation:
\beq
\label{TQZ}
{\bf T}(u){\bf Q}(u)=
\sin^{N}2\pi (u -\ell \gamma )
{\bf Q}(u +\gamma )
+\sin^{N}2\pi (u +\ell \gamma )
{\bf Q}(u -\gamma )\,.
\eeq

It is worth noting that the limit from (\ref{N6}) gives
an additional operator factor which is omitted in (\ref{Z4})
for it belongs to the center of the algebra generated by
the shifts $e^{2\gamma \p_{x_j}}$ and the
operator $\sum_{i=1}^{N}\p_{x_i}$. This factor has the form
$$
\prod_{i=1}^{N}
\frac{\Gamma_{\tilde q^2}(-2u\!+\!2\ell \gamma \!+\!1)
\Gamma_{\tilde q^2}(x_{i+1}\!-\!x_i \!+\!2\ell \gamma \!+\!1)}
{\Gamma_{\tilde q^2}(x_{i+1}\!-\!x_i \!-\!2u \!+\!1)}\,,
$$
where $\tilde q =\exp \left (-\frac{\pi i}{2\gamma}\right )$
is the modular partner of the $q=\exp (2\pi i \gamma )$
(cf.\,\cite{modular}). It seems so that the non-degenerate
Sklyanin algebra knows about both halves of the modular double.

\section{Concluding remarks}

In this work, using Baxter's method of vacuum vectors,
we have constructed a family of commuting difference
operators in many variables. This family provides the
$Q$-operator for the quantum integrable spin chain
of $XYZ$-type in the functional realization of the
algebra of observables.

The form of the obtained $Q$-operators
suggests a simple derivation of local integrals
of motion for the $XYZ$ chain in infinite-dimensional
representations. For instance, the natural candidates
for local Hamiltonians are
$({\bf Q}^{\pm}(\mp \ell \eta ))^{-1}
d/d\lambda \Bigl |_{\lambda = \mp \ell \eta}
{\bf Q}^{\pm}(\lambda )$.
Their explicit form is, however, not very illuminating.
The hermiticity properties
of these operators are to be further investigated.
We leave this for a future work.

Another important point remaining out of the scope
of this paper is the limit to the $XXX$-case
($\gamma \to 0$ in (\ref{Z1})). In this limit
one should be able to link our results with the
ones of the papers
\cite{derk,pronko,kss}, where the $Q$-operator
for the $XXX$ spin chain and allied models
with rational $R$-matrix was
realized as an integral operator.
However, it is
hard to compare them directly since the limit is
singular and, therefore, very
delicate. Moreover, as it should be clear
from the above, the rational limit is certainly not unique.
Different rational limits of the trigonometric
${\sf L}$-operator lead to ${\sf L}$-operators
for the $XXX$ spin chains, where the $sl(2)$-algebra
is realized in different ways. A detailed study of
rational limits seems to be interesting since in this case,
closely related to the representation theory of
$SL(2, \RR )$ or $SU(1,1)$,
the result must admit a meaningful translation into the
language of representation theory and harmonic analysis.
This may give a hint of how an elliptic extension of the
latter might look like.

The construction of the pre-$Q$-operators,
${\bf Q}_{L}(\lambda )$ and ${\bf Q}_{R}(\lambda )$, goes through
also for inhomogeneous chains. The problem, however,
is the normalization: in general there is no point
$\lambda_{*}$ such that
${\bf Q}_{L}(\lambda_{*})$ or
${\bf Q}_{R}(\lambda_{*})$ could be explicitly inverted.
Generally speaking, the kernel of the $Q$-operator
can not be represented as a product of local factors, so
local integrals of motion for inhomogeneous chains probably
do not exist. In other respects,
the pair of pre-$Q$-operators is a good substitute
for the $Q$-operator. Their kernels are factorized though
the operators do not commute. Presumably, they
still obey the exchange
commutation relations (\ref{Q1}). Their proof given
in Sec.\,6 should basically
go through in the inhomogeneous case,
too, since the key element of the proof is local in nature.

At last, we note that the technique
developed in the present paper might help to
construct the fundamental $R$-matrix in infinite-dimensional
representations of the Sklyanin algebra.

\section*{Acknowledgments}
I am grateful to E.\-An\-to\-nov, S.\-Khar\-chev,
S.\-Kho\-rosh\-kin,
I.\-Kri\-che\-ver, A.\-Le\-vin,
and T.\-Ta\-ke\-be for useful discussions, and to
V.\-Spi\-ri\-do\-nov for bringing the paper \cite{FrTur}
to my attention.
This work was supported in part by RFBR grant 98-01-00344.

\section*{Appendix A}
%\addcontentsline{toc}{section}{Appendix}
\def\theequation{A\arabic{equation}}
\setcounter{equation}{0}

\subsubsection*{Theta-functions}

We use the following definition of the
Jacobi $\theta$-functions:
\beq
\begin{array}{l}
\theta _1(z|\tau)=-\displaystyle{\sum _{k\in \z}}
\exp \left (
\pi i \tau (k+\frac{1}{2})^2 +2\pi i
(z+\frac{1}{2})(k+\frac{1}{2})\right ),
\\ \\
\theta _2(z|\tau)=\displaystyle{\sum _{k\in \z}}
\exp \left (
\pi i \tau (k+\frac{1}{2})^2 +2\pi i
z(k+\frac{1}{2})\right ),
\\ \\
\theta _3(z|\tau)=\displaystyle{\sum _{k\in \z}}
\exp \left (
\pi i \tau k^2 +2\pi i
zk \right ),
\\ \\
\theta _4(z|\tau)=\displaystyle{\sum _{k\in \z}}
\exp \left (
\pi i \tau k^2 +2\pi i
(z+\frac{1}{2})k\right ).
\end{array}
\label{theta}
\eeq
They also can be represented as infinite products.
The infinite product representation for the $\theta_1(z|\tau)$
reads:
\beq
\label{infprod}
\theta_1(z|\tau)=i\,\mbox{exp}\, \Bigl (
\frac{i\pi \tau}{4}-i\pi z\Bigr )
\prod_{k=1}^{\infty}
\Bigl ( 1-e^{2\pi i k\tau }\Bigr )
\Bigl ( 1-e^{2\pi i ((k-1)\tau +z)}\Bigr )
\Bigl ( 1-e^{2\pi i (k\tau -z)}\Bigr )
\eeq
Throughout the paper we write
$\theta _a(x|\tau)=\theta _a(x)$.
The transformation properties for shifts by the periods are:
\beq
\label{periods}
\theta_a (x\pm 1)=(-1)^{\delta _{a,1}+\delta _{a,2}}
\theta_a (x)\,,
\;\;\;\;\;
\theta_a (x\pm \tau )=(-1)^{\delta _{a,1}+\delta _{a,4}}
e^{-\pi i \tau \mp 2\pi i x}
\theta_a (x)\,.
\eeq
Under the modular transformation $\tau \to -1/\tau$
the $\theta$-functions behave as follows:
\beq
\label{mod}
\begin{array}{l}
\theta_{1}(z|\tau )=i \,\sqrt{i/\tau} \,
e^{-\pi i z^2/\tau }\theta_{1}(z/\tau |-1/\tau )\,,
\\ \\
\theta_{2}(z|\tau )= \sqrt{i/\tau} \,
e^{-\pi i z^2/\tau }\theta_{4}(z/\tau |-1/\tau )\,,
\\ \\
\theta_{3}(z|\tau )= \sqrt{i/\tau}\,
e^{-\pi i z^2/\tau }\theta_{3}(z/\tau |-1/\tau )\,,
\\ \\
\theta_{4}(z|\tau )= \sqrt{i/\tau} \,
e^{-\pi i z^2/\tau }\theta_{2}(z/\tau |-1/\tau )\,.
\end{array}
\eeq

By $\Theta_n$ we denote the space of $\theta$-functions
of order $n$, i.e., entire functions
$F(x)$, $x\in \CC$, such that
\beq
F(x+1)=F(x)\,,
\;\;\;\;\;\;
F(x+\tau)=(-1)^n e^{-\pi i n\tau -2\pi i nx}F(x)\,.
\label{8}
\eeq
It is easy to see that $\mbox{dim} \,\Theta_n =n$.
Let $F(x)\in \Theta_n$, then $F(x)$ has a multiplicative
representation of the form
$F(x)=c\prod _{i=1}^{n}\theta_1(x-x_i)$,
$\sum _{i=1}^{n}x_i =0$,
where $c$ is a constant. Imposing, in addition to (\ref{8}),
the condition $F(-x)=F(x)$, we define the space
$\Theta_{n}^{+}\subset \Theta_{n}$ of {\it even}
$\theta$-functions of order $n$, which
plays the important role in representations
of the Sklyanin algebra. If $n$ is an even number,
then $\mbox{dim}\, \Theta_{n}^{+} =\frac{1}{2}n +1$.

\subsubsection*{Elliptic gamma-function}

Here we collect the main formulas on the elliptic
gamma-function \cite{R3,FV3}. We use the (slightly modified)
notation of \cite{FV3}.
The elliptic gamma-function
is defined by the double-infinite product
\beq
\label{gamma}
\Gamma(z|\tau, \tau ')=\prod_{k,k'=0}^{\infty}
\frac{1-e^{2\pi i ((k+1)\tau +(k'+1)\tau ' -z)}}
{1-e^{2\pi i (k\tau +k'\tau ' +z)}}
\eeq
A sufficient condition for the product to be convergent is
$\mbox{Im}\,\tau >0$, $\mbox{Im}\,\tau' >0$.
We need the following properties of the elliptic gamma-function:
\beq
\label{gamma1}
\Gamma (z+1|\tau , \tau ')=
\Gamma (z|\tau , \tau ')\,,
\eeq
\beq
\label{gamma2}
\Gamma (z+\tau |\tau , \tau ')=
-ie^{-\frac{\pi i \tau'}{6}}\eta_{D}^{-1}(\tau')e^{\pi i z}
\theta_{1}(z|\tau')
\Gamma (z|\tau , \tau ')\,,
\eeq
\beq
\label{gamma3}
\Gamma (z+\tau' |\tau , \tau ')=
-ie^{-\frac{\pi i \tau}{6}}\eta_{D}^{-1}(\tau)e^{\pi i z}
\theta_{1}(z|\tau)
\Gamma (z|\tau , \tau ')\,,
\eeq
where
$$
\eta_{D}(\tau)=e^{\frac{\pi i \tau}{12}}
\prod_{k=1}^{\infty}\Bigl ( 1-e^{2\pi i k\tau}\Bigr )
$$
is the Dedekind function. Another useful property is
\beq
\label{gamma5}
\Gamma (z|\tau , \tau ')
\Gamma (\tau ' -z|\tau , \tau ')=
\frac{ie^{\pi i \tau '/6}\eta_{D}(\tau ')}{e^{\pi i z}
\theta_{1}(z|\tau ')}\,.
\eeq

Under the modular transformation $\tau \to -1/\tau$ the
elliptic gamma-function behaves as follows \cite{FV3}:
\beq
\label{modg}
\Gamma (z|\tau , \tau ')=e^{i \pi P(z)} \,
\frac{\Gamma (z/\tau \,|-1/\tau , \tau '/\tau )}
{\Gamma ((z-\tau )/\tau '\,|-\tau /\tau ', -1/\tau ')}\,,
\eeq
where
\beq
\label{polP}
\begin{array}{lll}
P(z)&=& \displaystyle{-\frac{1}{3\tau \tau '}\,z^3
+\frac{\tau +\tau ' -1}{2\tau \tau '}\,z^2
-\frac{\tau^2 +\tau '^2 +3\tau \tau ' -3\tau -3\tau ' +1}
{6\tau \tau '}\, z\,-} \\ &&\\
&-&\displaystyle{\frac{(\tau +\tau ' -1)(\tau +\tau ' -\tau \tau ')}
{12 \tau \tau '}}\,.
\end{array}
\eeq

Let us list the most frequently used
formulas for
$\Gamma (z) \equiv \Gamma (z|\tau , 2\eta )$.
Using (\ref{gamma3}) several times, we obtain:
\beq
\label{gam6}
\frac{\Gamma (x+2k\eta )}{\Gamma (x)}=
e^{\pi i \eta k^2}R^{-k}e^{\pi i kx}
\prod_{j=0}^{k-1}\theta_1 (x+2j\eta )\,,
\eeq
\beq
\label{gam7}
\frac{\Gamma (x-2k\eta )}{\Gamma (x)}=
(-1)^k e^{\pi i \eta k^2}R^{k}e^{-\pi i kx}
\prod_{j=0}^{k-1}
\Bigl ( \theta_1 (-x+2\eta +2j\eta )\Bigr )^{-1}\,,
\eeq
where $R=ie^{\pi i (\eta +\tau /6)}\eta_{D}(\tau )$.

As is seen from (\ref{gamma}), the function
$\Gamma (z|\tau , 2\eta )$ has zeros at the points
$z= 2(k+1)\eta + (m+1)\tau +n$,
and simple poles at the points
$z=-2k\eta -m\tau +n$,
where $k,m$ run over non-negative integers
and $n$ over all integers.
The residues of the elliptic gamma-function at the poles
at $z=-2k\eta$, $k=0,1,2, \ldots$ are:
\beq
\label{residue}
\mbox{res}\,\Bigl |_{z=-2k\eta}
\Gamma (z)=(-1)^k
e^{\pi i \eta k^2}R^{k}r_0
\prod_{j=1}^{k}
\Bigl ( \theta_1 (2j\eta )\Bigr )^{-1}\,,
\eeq
where
$$
r_0=
\mbox{res}\,\Bigl |_{z=0}
\Gamma (z)=
-\frac{e^{\pi i(\tau +2\eta )/12}}{2\pi i
\eta_D(\tau)\eta_D(2\eta )}\,.
$$
The function $\Phi_{\ell}(z)$
(\ref{Phi1}) can be extended to
arbitrary complex values of $\ell$ as follows:
\beq
\label{Phig}
\Phi_{\ell}(z)=e^{-\pi i \eta (2\ell +1)}R^{2\ell +1}
e^{-\pi i (2\ell +1)z}
\frac{\Gamma (z +2(\ell +1) \eta )}{\Gamma
(z -2\ell \eta )}\,.
\eeq

\subsubsection*{Elliptic hypergeometric series}

Here we follow \cite{FrTur}.
We define the elliptic Pochhammer symbol (the
shifted elliptic factorial) by
\beq
\label{Pohg} [x]_k \equiv [x] [x+1] \ldots [x+k-1]\,,
\eeq
where $[x]=\theta_1(2x\eta )$ (cf. (\ref{binom})).
By definition, the elliptic hypergeometric series is
\beq
\label{mh1}
{}_{r+1}W_{r}(\alpha_1 ; \alpha_4 , \alpha_5 , \ldots ,
\alpha_{r+1};z|\eta , \tau )=
\! \sum_{k=0}^{\infty}z^k
\frac{[\alpha_1 +2k][\alpha_1 ]_k}{[\alpha_1][k]!}
\prod_{m=1}^{r-2}
\frac{[\alpha_{m+3}]_k}{[\alpha_1 \!-\!\alpha_{m+3}\!+\!1]_k}\,.
\eeq
This is an elliptic analog of the very-well-poised basic
hypergeometric series \cite{GR}. The series is said to be
{\it balanced} if $z=1$ and
\beq
\label{mh2}
r-5+(r-3)\alpha_1 =2\sum_{m=1}^{r-2}\alpha_{m+3}\,.
\eeq
For a series $\sum_{k\geq 0}c_k$
of the form (\ref{mh1}), the balancing condition
(\ref{mh2}) means that the ratio
$c_{k+1}/c_k$ of the coefficients
is an elliptic function of $k$.
For balanced series (\ref{mh1}), we drop the
argument $z=1$ and the parameters $\eta , \tau$ writing
it simply as
${}_{r+1}W_{r}(\alpha_1 ; \alpha_4 , \ldots ,
\alpha_{r+1})$.
For instance,
\beq
\label{mh3}
{}_{10}W_{9}(\alpha_1 ; \alpha_4 , \alpha_5 , \ldots ,
\alpha_{10})=
\! \sum_{k=0}^{\infty}
\frac{[\alpha_1 +2k][\alpha_1 ]_k}{[\alpha_1][k]!}
\prod_{m=1}^{7}
\frac{[\alpha_{m+3}]_k}{[\alpha_1 \!-\!\alpha_{m+3}\!+\!1]_k}\,.
\eeq
The series is called
{\it terminating} if at least one of the parameters
$\alpha_4 , \ldots , \alpha_{r+1}$ is equal to a negative
integer number. In this case there is no problem of convergence.
If, say $\alpha_{r+1}=-n$, then the series terminates at $k=n$.
The terminating balanced
hypergeometric series were shown \cite{FrTur} to possess nice
modular properties.
That is why they were called {\it modular hypergeometric series}.

The modular hypergeometric series obey a number of
spectacular identities.
One of them is the elliptic analog
of the Bailey transformation for the ${}_{10}W_{9}$ \cite{GR}:
\beq
\label{Bailey}
\begin{array}{ll}
&{}_{10}W_{9}(\alpha_1 ; \alpha_4 , \ldots , \alpha_{10})
\\
=&{}_{10}W_{9}(\beta_1 ; \beta_4 , \ldots , \beta_{10})\,
\displaystyle{\frac{[\alpha_1 \!+\!1]_n[\beta_1 \!-
\!\beta_7 \!+\!1]_n [\beta_1 \!-\!
\beta_8 \!+\!1]_n [\alpha_1 \!-
\!\alpha_7 \! -\!\alpha_8 \!+\!1]_n}{[\beta_1 \!+\!1]_n[\alpha_1 \!-
\!\alpha_7 \!+\!1]_n [\alpha_1 \!-\!
\alpha_8 \!+\!1]_n [\beta_1 \!-
\!\beta_7 \! -\!\beta_8 \!+\!1]_n}}\,,
\end{array}
\eeq
where $\alpha_j$ obey the balancing condition $2+3\alpha_1 =
\sum_{m=1}^{7}\alpha_{m+3}$,
$\alpha_{10}=-n$, and the
parameters $\beta_j$ are as follows:
\beq
\label{beta}
\begin{array}{ll}
\beta_1 = 2\alpha_1 +1 -\alpha_4 -\alpha_5 -\alpha_6\,,&\\
\beta_i =\beta_1 -\alpha_1 +\alpha_i\,, &i=4,5,6, \\
\beta_i=\alpha_i \,, &i=7, \ldots , 10\,.
\end{array}
\eeq
This identity was proved in \cite{FrTur}.

\section*{Appendix B}
\def\theequation{B\arabic{equation}}
\setcounter{equation}{0}

Here we present the proofs of the intertwining relations
(\ref{S3}) and the commutation relation (\ref{T5}).

For general complex values of $\ell$ the operator
${\bf w}_{\ell}$ is given by the series (\ref{winf}):
\beq
\label{winf1}
{\bf w}_{\ell}=
\frac{e^{2\pi i (2\ell +1)z}\Gamma (2z)}{\Gamma (2z
+2(2\ell +1)\eta )}
\normord
{}_{4}W_{3}\Bigl (-\frac{z}{\eta}\!-\!2\ell \!-\!1 ;\,
-\!2\ell \!-\!1 ; \,e^{-2\eta \p_{z}}\Bigl | \eta , \tau \Bigr )
\normord
e^{(2\ell +1)\eta \p_{z}}
\eeq
The proof of (\ref{S3}) is a straightforward computation.
Substituting (\ref{Sa}), (\ref{winf1}) into (\ref{S3}) and
comparing the coefficients of $e^{2(\ell +1 -k)\eta \p_z}$,
$k>0$, we get, after cancelling out
common factors:
$$
\begin{array}{lll}
&&\theta_1 (2(\ell \!+\!k\!+\!1)\eta )
\theta_{a\!+\!1} (2z\!+\!2(\ell \!+\!1)\eta )
\theta_1 (2z \!+\!4k\eta )
\theta_{1} (2z\!-\!2(\ell \!-\!k \!+\!1)\eta )\\
\\ && \\
&+&
\theta_1 (2(\ell \!-\!k\!+\!1)\eta )
\theta_{a\!+\!1} (-2z\!+\!2(\ell \!+\!1)\eta )
\theta_1 (2z \!+\!4k\eta )
\theta_{1} (2z\!+\!2(\ell \!+\!k \!+\!1)\eta )\\
\\ && \\
&=&
\theta_1 (2(\ell \!+\!k\!+\!1)\eta )
\theta_{a\!+\!1} (2z\!+\!2(2k \!-\!\ell \!-\!1)\eta )
\theta_1 (2z )
\theta_{1} (2z\!+\!2(\ell \!+\!k \!+\!1)\eta )\\
\\ && \\
&+&
\theta_1 (2(\ell \!-\!k\!+\!1)\eta )
\theta_{a\!+\!1} (-2z\!-\!2(2k \!+\!\ell \!+\!1)\eta )
\theta_1 (2z )
\theta_{1} (2z\!-\!2(\ell \!-\!k \!+\!1)\eta )\,.
\end{array}
$$
All the four terms, as functions of $2z$,
have the same monodromy properties. They
are $\theta$-functions of 3-d order. The linear space
of such functions is 3-dimensional, so it is enough to
verify the equality at 3 points, say $z=0$, $z=-2k\eta$,
and $z=(\ell -k +1)\eta$. The coefficients of the
highest shift $e^{2(\ell +1 )\eta \p_{z}}$ in both sides
of (\ref{S3}) are equal, too, so the intertwining
relation is proved.

To prove (\ref{T5}),
we use the representation (\ref{fact}). For any two adjacent sites
of the chain (say $1$ and $2$) we define the $c$-number matrix
$\stackrel{12}{M}=
\stackrel{12}{M}(\lambda ,\ell)=
V^{-1}(\lambda_{-}, z_1) V(\lambda_{+}, z_2)$ explicitly
given by
$$
\stackrel{12}{M}(\lambda ,\ell)\!=\!\frac{1}{\theta_{1}(2z_2)}\!
\left (\!
\begin{array}{lll}
\theta_{1}(z_1 \!+\!z_2 \!-\!2\ell \eta )
\theta_{1}(z_1 \!-\!z_2 \!+\!2\lambda ) &&
\!-\theta_{1}(z_1 \!-\!z_2 \!-\!2\ell \eta )
\theta_{1}(z_1 \!+\!z_2 \!+\!2\lambda )
\\&&\\
\theta_{1}(z_1 \!-\!z_2 \!+\!2\ell \eta )
\theta_{1}(z_1 \!+\!z_2 \!-\!2\lambda ) &&
\!-\theta_{1}(z_1 \!+\!z_2 \!+\!2\ell \eta )
\theta_{1}(z_1 \!-\!z_2 \!-\!2\lambda )
\end{array}\!
\right )
$$
Using the cyclic property of the trace, one can rewrite
(\ref{T3}) as follows:
$$
{\bf T}(\lambda )\!=\!\normord \mbox{tr} \left [
\stackrel{12}{M}
\left ( \begin{array}{cc}
e^{\eta \p_{2}} &0
\\ 0& e^{-\eta \p_{2}}
\end{array} \right )
\stackrel{23}{M}
\left ( \begin{array}{cc}
e^{\eta \p_{3}} &0\\ 0 &e^{-\eta \p_{3}}
\end{array} \right )
\, \ldots \,
\stackrel{N1}{M}
\left ( \begin{array}{cc}
e^{\eta \p_{1}} &0\\ 0 &e^{-\eta \p_{1}}
\end{array} \right ) \!
\right ] \!
\normord
$$
where $\p_i =\p_{z_i}$. Computing the trace, we obtain:
\beq
\label{T4}
{\bf T}(\lambda) = \sum_{i_1 ,\ldots , i_N =1,2}
\!\stackrel{12}{M}_{i_1 i_2}
\stackrel{23}{M}_{i_2 i_3}
\ldots
\stackrel{N1}{M}_{i_N i_1}
e^{(3-2i_1)\eta \p_{1}}
e^{(3-2i_2)\eta \p_{2}} \ldots
e^{(3-2i_N)\eta \p_{N}}\,.
\eeq
Recall that
$\varphi_{\ell}(z_i , z_{i+1})=
\Phi_{\ell}(z_i \!-\!z_{i+1})\Phi_{\ell}(z_i \!+\!z_{i+1})$.
The commutation relation (\ref{T5}),
$$
{\bf T}^{(\ell )}(\lambda)\,\prod_{i=1}^{N}\varphi_{\ell}(z_i , z_{i+1})
=\Bigl (\prod_{i=1}^{N}\varphi_{\ell}(z_i , z_{i+1})\Bigr )
{\bf T}^{(-\ell -1)}(\lambda)\,,
$$
follows from (\ref{T4}) and the local
commutation rule
$$
\begin{array}{lll}
&&e^{(3-2i_1)\eta \p_{1}\! +\!(3-2i_2)\eta \p_{2}}
\varphi_{\ell}(z_1 , z_{2}) \\ &&\\
&=&\displaystyle{\frac{\theta_1 (z_1+(-1)^{i_1 +i_2}z_2 \!-
\!2(-1)^{i_1}(\ell +1)\eta )}{\theta_1 (z_1+(-1)^{i_1 +i_2}z_2 \!+
\!2(-1)^{i_1}\ell \eta )}}\,
\varphi_{\ell}(z_1 , z_{2})
e^{(3-2i_1)\eta \p_{1}\! +\!(3-2i_2)\eta \p_{2}}\,,
\;\;i_1 , i_2 =1,2\,,
\end{array}
$$
which is easy to verify.
Writing out the entries of the matrix $M$ in the form
\beq
\label{T6}
\stackrel{12}{M}_{ik}=\theta_1 \Bigl (z_1 +(-1)^{i+k}z_2
+2(-1)^i \ell \eta \Bigr )
C_{ik} \,,
\;\;\;\;\;\;i,k =1,2\,,
\eeq
(the second factor stands for something that does not
depend on $\ell$), and subsequently moving the factors
$\varphi_{\ell}(z_i , z_{i+1})$
from the right to the left through
each term of (\ref{T4}), we arrive at (\ref{T5}).

\end{document}